\documentclass[10pt]{amsart}
\usepackage{amssymb,amsmath,amsthm}
\usepackage{mathrsfs,dsfont,a4wide}
\theoremstyle{plain}
\newtheorem{theorem}{Theorem}[section]
\newtheorem{lemma}[theorem]{Lemma}
\newtheorem{proposition}[theorem]{Proposition}

\newtheorem{remark}[theorem]{Remark}
\newtheorem{definition}[theorem]{Definition}
\theoremstyle{definition}
\theoremstyle{remark}
\numberwithin{equation}{section}

\newcommand{\Tr}{Tr}
\newcommand{\SymN}{{\mathcal{S}^{N\times N}}}

\newcommand{\as}{{\mathcal A}}
\newcommand{\hs}{{\mathcal H}}
\newcommand{\ks}{{\mathcal K}}

\newcommand{\fs}{{\mathcal F}}

\newcommand{\bs}{{\mathcal B}}

\newcommand{\Es}{{\mathcal E}}

\newcommand{\rs}{{\mathcal R}}
\newcommand{\ws}{{\mathcal W}}

\newcommand{\R}{{\mathbb R}}
\newcommand{\N}{{\mathbb N}}

\newcommand{\tint}[1]{{\rm int}(#1)}

\newcommand{\diam}[1]{{\rm diam}(#1)}

\newcommand{\Om}{\Omega}
\newcommand{\Omb}{\overline{\Omega}}

\newcommand{\hn}{\hs^{N-1}}
\newcommand{\hu}{\hs^1}

\newcommand{\eps}{\varepsilon}

\newcommand{\Sg}[2]{S^{#1}(#2)}

\newcommand{\Div}{{\rm div}}
\title
[Crack initiation in elastic bodies]
{Crack initiation in elastic bodies}
\author[A. Chambolle]
{Antonin Chambolle}
\address[Antonin Chambolle]{
CMAP, Ecole Polytechnique,
91128 Palaiseau Cedex, France }
\email[A. Chambolle]{antonin.chambolle@polytechnique.fr}
\author[A. Giacomini]
{Alessandro Giacomini}
\address[Alessandro Giacomini]{Dipartimento di Matematica, Facolt\`a di Ingegneria, Universit\`a degli Studi di Brescia, Via Valotti 9, 25133 Brescia, Italy}
\email[A. Giacomini]{alessandro.giacomini@ing.unibs.it}
\author[M. Ponsiglione]
{Marcello Ponsiglione}
\address[Marcello Ponsiglione]{Max Planck Institute for
Mathematics in the Sciences, Inselstrasse 22, D-04103 Leipzig, Germany}
\email[M. Ponsiglione]{marcello.ponsiglione@mis.mpg.de}
\begin{document}
\vskip .2truecm
\begin{abstract}
\small{
In this paper we study the crack initiation in a hyper-elastic body
governed by a Griffith's type energy. 
We prove that, during a load process through a time dependent boundary datum of the type $t \to t g(x)$ and in absence of strong singularities (this is the case of homogeneous isotropic materials) the crack initiation is brutal, i.e., 
a big crack appears after a positive time $t_i>0$. On the contrary, in presence of a point 
$x$ of strong singularity, a crack will depart from $x$ at the initial time of
loading and with zero velocity. We prove these facts (largely expected by the experts
of material science) for admissible cracks belonging to the large class of closed one dimensional sets with a finite number of connected components.
\par
The main tool we employ to address the problem is a local minimality result for the functional
$$
\Es(u,\Gamma):=\int_\Om f(x,\nabla v)\,dx+k\hu(\Gamma),
$$
where $\Omega \subseteq \R^2$, $k>0$ and $f$ is a suitable Carath\'eodory function. We prove that if the uncracked configuration $u$ of $\Om$ relative to a boundary displacement $\psi$ has uniformly weak singularities, then configurations $(u_\Gamma,\Gamma)$ with 
$\hu(\Gamma)$ small enough are such that $\Es(u,\emptyset)<\Es(u_\Gamma,\Gamma)$.
}
\end{abstract}
\maketitle
{\small \tableofcontents}
 
\section{Introduction}
Griffith's criterion for crack propagation in hyper-elastic bodies asserts that, during a load process, a crack $\Gamma$ can grow only if the energy dissipated  to enlarge the crack, which is basically assumed to be proportional to the area of the cracked surface, is  balanced by the corresponding release of bulk energy. 
According to {\it Griffith's theory}, if $\Om$ represents a two dimensional hyper-elastic body, $\psi$ is a boundary datum and 
$\Gamma$ is a curve in $\Om$ parametrized by arc length,
then the crack $\Gamma(l_0)$ is in equilibrium if 
\begin{equation}
\label{griffitheq0}
k(l_0):=\limsup_{l \to 0^+} \frac{W(u(l_0))-W(u(l_0+l))}{l} \le k,
\end{equation}
where $u(l_0)$ and $u(l_0+l)$ are the displacements associated to 
$\psi$ and to the cracks $\Gamma(l_0)$ and $\Gamma(l_0+l)$ respectively, $W$ is the {\it bulk energy functional} and $k$ is the {\it toughness} of the material. A {\it quasistatic crack evolution} is determined by an increasing function $t \to l(t)$ satisfying the {\it Griffith's criterion} for crack propagation, which asserts that for every $t$ we have $k(l(t)) \le k$ and
$$
(k-k(l(t)))\dot{l}(t)=0,
$$
i.e., $\Gamma(l(t))$ propagates if and only if \eqref{griffitheq0} holds with equality.
\par
The aim of this paper is to discuss this criterion in the case of {\it crack initiation},
i.e., when there is not a pre-existing crack in the body ($l_0=0$).  
A fundamental role in the problem is played by the singularities of the body, namely the behavior of the elastic energy concentration of the deformation. Experiments show that small cracks usually appear near sufficiently strong singular points of the the body, whose position are essentially determined by its inhomogeneities. If the singularities of the body are sufficiently weak (for instance this is the case of homogeneous isotropic materials), a lot of results in the literature of material science show that the derivative in \eqref{griffitheq0} for $l_0=0$ is equal to zero. The conclusion is that Griffith's criterion is  not adequate to predict crack initiation (and, as a consequence, a crack evolution originating from an uncracked configuration). These results require that the path of the crack is sufficiently regular (a line or a smooth curve). In this paper we prove that the same conclusion holds in the class of all one dimensional closed sets with a finite number of connected components. More precisely we prove that the limit in \eqref{griffitheq0} is zero if $\Gamma(l)$ is any family of closed sets with length less than $l$
and with at most $m$ connected components, with $m$ independent of $l$. In particular we do not prescribe the path nor the shape of the cracks. 
\par
Although it is more general,
our study is in part motivated by the variational model for quasistatic
crack propagation proposed by Francfort and Marigo in~\cite{FM}.
The main features of this model
are that the path of the increasing crack $\Gamma(t)$
is not preassigned, the class of admissible cracks is given by all sets with finite length, and the growth is not assumed to be progressive, namely the length of the crack is not assumed to be continuous in time.
The classical Griffith's equilibrium condition
for the configuration $(u(t),\Gamma(t))$ is replaced
by a {\it unilateral minimality} property and an {\it energy balance} condition. 
The unilateral minimality property states that, during the crack evolution, the total energy is minimal among all configurations with larger cracks, namely
\begin{equation} 
\label{ump}
W (u(t)) + k\hu(\Gamma(t)) \le W(v) + k\hu(H),
\end{equation}
for every crack $H$ larger than $\Gamma(t)$ and for every deformation $v$ admissible for the boundary datum $\psi$ and for $H$. (Here $\hu$ --- the 1-dimensional Hausdorff measure --- is a suitable generalization of the length).
The energy balance condition states that the energy of the system evolves in relation with the power of external loads in such a way that no dissipation occurs.
The authors claim that their model improves the understanding of the crack initiation with respect to the classical Griffith's criterion: as a matter of fact,
in contrast with Griffith's model, it
admits {\it brutal crack initiation}, i.e., evolutions $\Gamma(t)$ of the type
$$
\Gamma(t)=\emptyset 
\qquad\text{ for every }t\le t_i 
$$
and 
$$
\inf_{t>t_i}\hs^1(\Gamma(t))>0,
$$  
where $t_i$ is referred to as {\it time initiation} of the crack.  
In this paper we prove that, within the class of cracks which are closed and with at most $m$ connected components, crack initiation is always brutal whenever the elastic displacement presents sufficiently weak singularities. More precisely we show that \eqref{ump} can be violated only by cracks whose length is greater than a critical quantity $l^*$, depending on the boundary datum and on the physical properties of the material. On the contrary, in presence of a point $x$ of strong singularity the crack initiation is progressive: a crack departs from $x$ at the initial time of loading and  with zero velocity. 
These facts were proved by Francfort and Marigo in \cite[Proposition 4.19]{FM}, under the assumption that the path of the crack is given a priori by a finite number of fixed curves which can be parametrized by arc length. This is not the case in our larger class of admissible cracks.
\par
The main tool we employ to address the problem of crack initiation
%% in the theories of Griffith and Francfort-Marigo
is a local minimality result for the
% Mumford-Shah [EXPLAIN] type 
functional
\begin{equation}
\label{mainfunctional}
\int_\Om f(x,\nabla v)\,dx+k\hu(\Gamma),
\end{equation}
where $\Om$ is a bounded Lipschitz open set in $\R^2$, $k>0$, 
and $f: \Om \times \R^2 \to \R$ is
a Carath\'eodory function strictly convex and $C^1$ in the second variable,
satisfying standard $p$-growth estimates with $p>1$, and such that
$f(x,0)=0$. The functional~\eqref{mainfunctional} is a variant of a
functional which has first appeared in the theory of image segmentation,
in a celebrated paper by Mumford and Shah~\cite{MumSha}.
The set $\Gamma$ belongs to the class 
\begin{equation}
\label{defkmintr2}
\ks_m(\Omb):=\{\Gamma \subseteq \Omb\,:\, \Gamma \text{ has at most $m$ connected components and }\hu(\Gamma)<+\infty\}
\end{equation}
and the function $v$ belongs to the Sobolev space
$W^{1,p}(\Om \setminus \Gamma)$
and satisfies the boundary condition
\begin{equation}
\label{bdrycond}
v=\psi
\qquad
\text{on }\partial_D \Om \setminus \Gamma,
\end{equation}
where $\partial_D\Om \subseteq \partial \Om$ is open in the relative topology, and
$\psi$ is (the trace of a function) in $W^{1,p}(\Om) \cap L^\infty(\Om)$.
\par
Let $u_\Gamma$ be a minimum energy displacement relative to $\psi$ and $\Gamma$, i.e., let $u_\Gamma$ be a minimizer for
\begin{equation}
\label{elasticsol0bis}
\min \left\{ \int_\Om f(x,\nabla v)\,dx\,:\, u \in W^{1,p}(\Om \setminus \Gamma), v=\psi \text{ on }\partial_D \Om \right\}.
\end{equation}
We denote by $u$ the {\it elastic configuration} of $\Om$ relative to the boundary datum $\psi$, i.e., a solution of \eqref{elasticsol0bis} with $\Gamma=\emptyset$, and we assume that $u$ admits {\it uniformly weak singularities} in $\Omb$, i.e., 
\begin{equation}
\label{weaksingbis}
\|\nabla u\|^p_{L^p(B_r \cap \Om)} \le Cr^\alpha
\end{equation}
for some constants $\alpha>1$ and $C>0$ and for every ball $B_r$ with radius $r$.
Condition \eqref{weaksingbis} means that  the bulk energy of the elastic configuration $u$ in a ball $B_r(x)$ is negligible with respect to the length of $\partial B_r(x)$ as $r$ goes to zero, uniformly in $x \in \Omb$.
\par
Our main result is the following Theorem, which establishes that under the previous assumptions small cracks are not energetically convenient for the functional \eqref{mainfunctional}.

\begin{theorem}
\label{mainthmbis}
Assume that $u$ admits only uniformly weak singularities in $\Omb$. Then there exists a critical length $l^*$ depending only on $\Om$, $f$, $k$, $\psi$ and $m$ such that for all $\Gamma \in \ks_m(\Omb)$ with $\hu(\Gamma)<l^*$ we have
\begin{equation}
\label{stima0}
\int_\Om f(x,\nabla u)\,dx <\int_\Om f(x,\nabla u_\Gamma)\,dx+ k \hu(\Gamma).
\end{equation}
\end{theorem}

We observe (see Remark~\ref{l1topology}) that this statement is equivalent to the
local minimality of $u$ in~\eqref{mainfunctional}, in the $L^1$
topology.

Let us briefly comment the assumption about the singularities of $u$
in  Theorem \ref{mainthmbis}. The minimality result  is false if the elastic solution $u$ has {\it strong singularities}, namely if there exists $x \in \Omb$ such that
\begin{equation}
\label{ess0}
\limsup_{r \to 0} \frac{1}{r}\int_{B_r(x) \cap \Omb}|\nabla u|^p\,dx=+\infty.
\end{equation} 
In fact condition \eqref{ess0} ensures that it is energetically convenient to create a small crack $\Gamma:= \partial B_r(x)$ around $x$: the surface energy needed to create such a crack is proportional to $r$, while the corresponding release of bulk energy 
is by \eqref{ess0} bigger than $r$ if $r$ is small enough.
\par 
The critical case when the right hand-side of \eqref{ess0} is a constant $0<C<\infty$ corresponds to the  singularity appearing around the tip of the crack (see \cite{GRI}). In this case the celebrated Irwin's formula states that the release of bulk energy per unit length along rectilinear increments of the crack is equal to the so called {\it mode III stress intensity factor} $K_\textit{III}$, which is proportional
to $C$. In our class of cracks $\ks_m(\Omb)$ we have that if $C$ is small enough, 
then the release of bulk energy per unit length is less than $k$,
and therefore our minimality result still holds, while it is false if $C$ is too large. We can not fill the gap, and therefore we do not achieve a sharp Irwin type formula  in our  class of cracks.
\par
In order to prove Theorem \ref{mainthmbis} we have to compare the asymptotic behavior 
of the release of bulk energy
\begin{equation}
\label{release}
\int_\Om [f(x,\nabla u)-f(x,\nabla u_{\Gamma})]\,dx
\end{equation}
with $\hs^1(\Gamma)$ when $\hs^1(\Gamma)\to 0$. In the literature there are many results in this direction considering particular sequences of infinitesimal cracks $\Gamma_n$, for instance when $\Gamma_n$ is the rescaled version of a fixed smooth curve $\Gamma$. 
An intuitive strategy to estimate \eqref{release} is to compute how much energy 
is required in order to make $u_\Gamma$  a good competitor for the minimum problem   \eqref{elasticsol0bis} without cracks, namely how much energy is required to heal the crack $\Gamma$. 
\par
This seems difficult for a generic crack in $\ks_m(\Omb)$. So our strategy is to operate on the stress $\sigma:=\partial f(x,\nabla u)$ of the elastic solution. More precisely we prove the following key estimate (see \eqref{mainest}) 
\begin{equation}
\label{mainest0}
\int_\Om [f(x,\nabla u)-f(x,\nabla u_\Gamma)]\,dx \le 
\int_\Om [\tau -\sigma] \cdot [\partial f^*(x, \tau)-\partial f^*(x,\sigma)] \, dx,
\end{equation}
for all vector fields $\tau\in L^q(\Om;\R^N)$ ($q=p':=\frac{p}{p-1}$) such that
\begin{equation}
\label{divgammaintr}
\int_\Om \tau \cdot \nabla v=0 \text{ for all }v \in \as(\Gamma).
\end{equation}
Here $f^*$ is the convex conjugate of $f$, and $\as(\Gamma)$ is defined as
$$
\as(\Gamma):= \{v \in W^{1,p}(\Om \setminus \Gamma)\,:\, v=0 \text{ on }\partial_D \Om \setminus \Gamma \}.
$$
If $\Gamma$ is sufficiently regular, condition \eqref{divgammaintr} means that 
the vector field $\tau$ has zero divergence outside $\Gamma$, and $\tau(x)$ is tangent to $\Gamma$ for every $x \in \Gamma$.
\par
The proof of Theorem \ref{mainthmbis} relies on the construction of a vector field $\tau$ 
satisfying \eqref{divgammaintr} and such that
\begin{equation}
\label{desiredineq}
\int_\Om [\tau -\sigma] \cdot [\partial f^*(x, \tau)-\partial f^*(x,\sigma)]\,dx
<k\hu(\Gamma).
\end{equation}
We construct $\tau$ modifying the stress $\sigma$ which is divergence-free 
(as a consequence of Euler equation of problem \eqref{elasticsol0bis}), 
but  not tangent to $\Gamma$. First, we consider a neighborhood $U$ of $\Gamma$, and a cut-off function $\varphi$ such that $\varphi=0$ on $U$. Then
$\varphi \sigma$ is null near $\Gamma$, and in particular it is tangent to $\Gamma$.
Then we construct a vector field $\eta$ in such a way that $\eta=0$ on $U$ and
$$
{\rm div}\, \eta=-{\rm div}(\varphi \sigma).
$$
We get that $\tau=\varphi \sigma+\eta$ is an admissible vector fields for inequality \eqref{mainest0}. Using the fact that $\Gamma \in \ks_m(\Omb)$ and that $u$ has uniformly weak singularities, it is possible to choose $U$, $\varphi$ and $\eta$ in such a way that inequality \eqref{desiredineq} holds.
\par
It turns out that the constraint \eqref{divgammaintr} can be handled in an easier way than the constraint of being a gradient, and this is the reason why we work with the stress $\sigma$ instead of the strain $\nabla u$.
\par
Estimate \eqref{mainest0} actually holds true in any dimension, and it turns out that our arguments work in any dimension provided that the crack $\Gamma$ belongs to the class
\begin{equation}
\label{ksCintr}
\ks^C(\Omb):=\{\Gamma \subseteq \Omb\,:\, \Gamma \text{ is closed and } \diam{\Gamma} \le C\hn(\Gamma)\},
\end{equation}
where $C$ is a fixed constant.
This is certainly true in dimension two for the cracks in $\ks_1(\Omb)$
(and in $\ks_m(\Omb)$ up to a localization argument).
However, it also shows that the local minimality result of Theorem
\ref{mainthmbis} remains valid in higher dimension,
within the class $\ks^C(\Omb)$ of cracks that are {\it not needle-like}:
see Remark~\ref{higherdim}. 
\par
The minimality result holds also in the case of planar linearized
elasticity, with a density of bulk
energy involving the symmetrized gradient. This is considered in
Section~\ref{planarelasticitysec} --- while the ``simpler'' case of
2D vectorial nonlinear elasticity is addressed in Remark~\ref{secvectorial}.
\par
Finally in Appendix \ref{SBVappendix} we extend the local minimality result of Theorem \ref{mainthmbis} to the larger class of all 1-dimensional rectifiable sets.
% the displacements of class $SBV$.
This seems to be the most general class: however,
there is a price to pay in order to handle such admissible
cracks.
First, we need to assume that
$\nabla u$ is bounded and regular up to the boundary.
In particular we are not able to treat the case in
which $u$ has weak singularities.
Then, the method we employ is based on the maximum principle, which allows
to estimate the local opening of a crack with the global energy in a
small ball surrounding the crack. It is therefore
strictly scalar, and bidimensional.
\smallskip
\par
The paper is organized as follows. In Section~\ref{mainestsec} we establish the main inequality \eqref{mainest0}. In Section~\ref{mainresultsec} we prove the local minimality result in dimension $2$ and its extensions to vector-valued displacements and 
to the $N$-dimensional case within the class of cracks given by \eqref{ksCintr}. Section~\ref{planarelasticitysec} addresses the case of planar two-dimensional elasticity.
The problem of crack initiation in quasistatic evolutions is addressed
in Section~\ref{inisec}, while the two-dimensional $SBV$-case without singularities is treated in  Appendix \ref{SBVappendix}. In Appendix \ref{app:PK} we show how to obtain some uniform Poincar\'e and Poincar\'e-Korn type inequalities,
used during the proofs of our main results in the construction of the competitor stress field $\tau$.

\section{The dual problem and the main estimate}
\label{mainestsec}
Let $\Om$ be a bounded connected Lipschitz open set in $\R^N$, let $\partial_D \Om \subseteq \partial \Om$ be open in the relative topology, and let $\partial_N \Om:=\partial \Om \setminus \partial_D \Om$. Let $f: \Om \times \R^N \to \R$ be a Carath\'eodory function such that
\begin{equation}
\label{fstruttura}
\xi \to f(x,\xi) \text{ is strictly convex and $C^1$ for a.e. }x \in \Om,
\end{equation}
\begin{equation}
\label{fzero}
f(x,0)=0
\qquad
\text{for a.e. }x \in \Om,
\end{equation}
and such that for a.e. $x \in \Om$ and for all $\xi \in \R^N$
\begin{equation}
\label{fcrescita}
\alpha |\xi|^p \le f(x,\xi) \le \beta(|\xi|^p+1),
\end{equation}
where $\alpha,\beta>0$ and $1<p<+\infty$.
\par
Given $\psi \in W^{1,p}(\Om) \cap L^\infty(\Om)$ and $\Gamma$ closed set contained in $\Omb$, let us consider the minimization problem
\begin{equation}
\label{elasticsol}
\min \left\{ \int_\Om f(x,\nabla u)\,dx\,:\, u \in W^{1,p}(\Om \setminus \Gamma), 
u=\psi \text{ on }\partial_D \Om \setminus \Gamma \right\}.
\end{equation}
In view of \eqref{fzero} and since $\psi \in W^{1,p}(\Om) \cap L^\infty(\Om)$, problem \eqref{elasticsol} is well posed.
\par
Let us denote by $u_\Gamma \in W^{1,p}(\Om \setminus \Gamma) \cap L^\infty(\Om \setminus \Gamma)$ a minimizer of \eqref{elasticsol}. Clearly $\nabla u_\Gamma$ is uniquely determined, while $u_\Gamma$ is determined up to a constant on each connected component of $\Om \setminus \Gamma$ which does not touch $\partial_D \Om$. 
\par
We denote by $u \in W^{1,p}(\Om) \cap L^\infty(\Om)$ the solution of \eqref{elasticsol} corresponding to $\Gamma=\emptyset$, and we refer to $u$ as the {\it elastic solution}. In the case $\Gamma$ is sufficiently regular,  the Euler Lagrange equation satisfied by $u_\Gamma$ is
\begin{equation}
\label{elequation}
\begin{cases}
\Div \partial_\xi f(x,\nabla u_\Gamma)=0      & \text{on }\Om \setminus \Gamma, \\
u=\psi      & \text{on }\partial_D \Om \setminus \Gamma, \\
\partial_\xi f(x,\nabla u_\Gamma) \cdot n=0 &\text{on }\partial_N \Om \cup \partial \Gamma,
\end{cases}
\end{equation}
where $n$ denotes the normal vector to $\partial_N \Om \cup \Gamma$.
In the sequel, we will write $\partial f(x,\xi)$ for $\partial_\xi f(x,\xi)$. 
\par
Let us set
\begin{equation}
\label{agamma}
\as(\Gamma):= \{v \in W^{1,p}(\Om \setminus \Gamma)\,:\, v=0 \text{ on }\partial_D \Om \setminus \Gamma \}.
\end{equation}
\par
Let us denote by $f^*$ the convex conjugate of $f$ with respect to the second variable
defined by
$$
f^*(x,\zeta):=\sup \{\zeta\cdot\xi-f(x,\xi)\,:\, \xi \in \R^N\}.
$$
We refer the reader to \cite{R} for the main properties of the conjugate function $f^*$. Notice that $f^*$ is of class $C^1$ since $f$ is strictly convex. The main result of this section is the following.

\begin{theorem}
\label{mainestthm}
Let $\Gamma$ be a closed subset of $\Omb$, and let $\sigma:=\partial f(x,\nabla u)$ be the stress associated to the elastic configuration $u$. Then we have
\begin{equation}
\label{mainest}
\int_\Om [f(x,\nabla u)-f(x,\nabla u_\Gamma)]\,dx \le \int_\Om 
[\tau -\sigma] \cdot [\partial f^*(x, \tau)-\partial f^*(x,\sigma)]\,dx
\end{equation}
for all $\tau \in L^q(\Om;\R^N)$ ($q=p':=\frac{p}{p-1}$) such that
\begin{equation}
\label{divgamma}
\int_\Om \tau \cdot \nabla v\,dx=0 \text{ for all }v \in \as(\Gamma).
\end{equation}
\end{theorem}

\begin{proof}
For all $\eta \in L^p(\Om;\R^N)$ let us set
\begin{equation}
\label{phifunctional}
\Phi(\eta):= \min_{w \in u+\as(\Gamma)} \int_\Om f(x,\nabla w+\eta)\,dx,
\end{equation}
where $\as(\Gamma)$ is defined in \eqref{agamma}. Then the convex conjugate of $\Phi$ defined on $L^q(\Om;\R^N)$ has the form
\begin{align}
\label{phiconj}
\Phi^*(\tau) &:= \sup_{w,\eta} \int_\Om [\tau \cdot \eta -f(x,\nabla w+\eta)] \,dx \\
\nonumber &=\sup_{w,\eta} \int_\Om [\tau \cdot (\eta+\nabla w) -f(x,\nabla w+\eta) -\tau \cdot \nabla w] \,dx \\
\nonumber &=\int_{\Om}[f^*(x,\tau)-\tau \cdot \nabla u]\,dx + \sup_{v \in \as(\Gamma)}\int_\Om \tau \cdot \nabla v\,dx.
\end{align}
We conclude that
\begin{equation}
\label{phiconj2}
\Phi^*(\tau)=
\begin{cases}
\displaystyle \int_\Om [f^*(x,\tau)-\tau \cdot \nabla u]\,dx&\text{if $\tau$ satisfies}\eqref{divgamma} \\
+\infty &\text{otherwise.}
\end{cases}
\end{equation}
\par\noindent
Notice that $\Phi(0)=\Phi^{**}(0)$ because $\Phi$ is weakly lower semicontinuous and 
$$
\Phi(0)=\int_\Om f(x,\nabla u_\Gamma)\,dx<+\infty.
$$ 
(In fact, $\Phi$ is locally finite, hence locally Lipschitz).
Therefore  we obtain
$$
- \int_\Om f(x,\nabla u_\Gamma)\,dx = - \Phi(0) =
- \Phi^{**}(0) = \min_{\tau} \Phi^*(\tau),
$$
so that by \eqref{phiconj2} we deduce
$$
-\int_\Om f(x,\nabla u_\Gamma)\,dx=\min_\tau
\left\{ \int_\Om [f^*(x,\tau)-\tau \cdot \nabla u]\,dx\,:\,  \tau \text{ satisfies \eqref{divgamma}} \right\}. 
$$
For all $\tau$ satisfying \eqref{divgamma},  we get
$$
\int_\Om [f(x,\nabla u)-f(x,\nabla u_\Gamma)]\,dx \le 
\int_\Om [f(x, \nabla u)+f^*(x,\tau)-\tau \cdot \nabla u]\,dx.
$$
Let $\sigma(x):= \partial f(x,\nabla u(x))$ be the stress of the elastic solution $u$. Since  for a.e. $x \in \Om$
$$ 
f^*(x,\tau(x)) \le f^*(x,\sigma(x))+ \partial f^*(x,\tau(x)) \cdot(\tau(x)-\sigma(x)),
$$
and
$$
f(x,\nabla u(x))+f^*(x,\sigma(x))=\nabla u(x) \cdot \sigma(x),
$$ 
and since $\nabla u(x)=\partial f^*(x,\sigma(x))$, we finally obtain our main estimate \eqref{mainest}, so that the proof is concluded.
\end{proof}

\begin{remark}
\label{explaindivgamma}
{\rm
If $\Gamma$ is sufficiently regular, condition \eqref{divgamma} means that $\tau$ has zero divergence outside $\Gamma$, and $\tau(x)$ is tangent to $\Gamma$ for every $x \in \partial \Gamma$.
}
\end{remark}

\section{The minimality result in anti-plane elasticity}
\label{mainresultsec}
In this section we prove that under some assumptions on the elastic configuration $u$, small cracks are not convenient for the total energy
\begin{equation}
\label{totalenergy}
\int_\Om f(x, \nabla v)\,dx+k\hu(\Gamma),
\end{equation}
where $f$ is a Carath\'eodory function satisfying conditions \eqref{fstruttura}, \eqref{fzero} and \eqref{fcrescita}, and $k>0$.
\par
Let us consider $\Om$ bounded connected Lipschitz open subset of $\R^2$, and let $\partial_D \Om \subseteq \partial \Om$ be open in the relative topology and such that $\partial_N \Om:=\partial \Om \setminus \partial_D \Om$ has a finite number of connected components.
\par
Let $m$ be a positive integer. The class of admissible cracks is given by
\begin{equation}
\label{defkm}
\ks_m(\Omb):=\{\Gamma \subseteq \Omb\,:\, \Gamma \text{ has at most $m$ connected components and }\hu(\Gamma)<+\infty\}.
\end{equation}
\par
Given $\psi \in W^{1,p}(\Om) \cap L^{\infty}(\Om)$ and $\Gamma \in \ks_m(\Omb)$, the displacement 
$$
u_\Gamma \in W^{1,p}(\Om \setminus \Gamma)
$$ 
associated to $\Gamma$ and $\psi$ is given by problem \eqref{elasticsol}. We denote with $u$ the solution of \eqref{elasticsol} relative to $\Gamma=\emptyset$, and we refer to $u$ as the {\it elastic solution}. 
\par
The basic assumption on the elastic configuration $u$ involves the behavior of the energy concentration of the stress. We require that $u$ has {\it uniformly weak singularities} in $\Omb$, in the sense of the following definition.

\begin{definition}
\label{weaksingdef}
{\rm
We say that $u \in W^{1,p}(\Om)$ has only {\it uniformly weak singularities} in $A$, open subset of $\Om$, if there exist constants $1<\alpha<2$ and $C>0$ such that for every $x \in A$ and for every $r$ small
\begin{equation}
\label{weaksing}
\int_{B_r(x) \cap \Om}|\nabla u|^p\,dx \le Cr^\alpha.
\end{equation}
}
\end{definition}

As mentioned in the Introduction, the condition of uniformly weak singularities means that 
the bulk energy of the elastic configuration $u$ in a ball $B_r(x)$ is an infinitesimal of higher order than the length of $\partial B_r(x)$ as $r$ goes to zero, uniformly with respect to $x \in A$.
\par
The main result of this section is the following theorem.

\begin{theorem}
\label{mainthm}
Let the elastic solution $u$ of problem \eqref{elasticsol} have only uniformly weak singularities in $\Om$ according to Definition \ref{weaksingdef}.  Then there exists a critical length $l^*>0$ depending only on $\Om$, $m$, $f$, $k$ and $\psi$ such that for all $\Gamma \in \ks_m(\Omb)$ with $\hu(\Gamma)<l^*$ we have
\begin{equation}\label{mmc}
\int_\Om f(x,\nabla u)\,dx< \int_\Om f(x,\nabla u_{\Gamma})\,dx+k\hu(\Gamma),
\end{equation}
where $u_{\Gamma}$ is a minimum of \eqref{elasticsol}.
\end{theorem}

\begin{remark}
\label{l1topology}
{\rm {\bf (Local minimality in the $L^1$-topology)}
The minimality condition \eqref{mmc} implies that
the elastic solution $u$ is a local minimum for the total energy \eqref{totalenergy}
with respect to the $L^1$-topology. More precisely, for every sequence $(\Gamma_h)_{h \in \N}$ in $\ks_m(\Omb)$ and for every $u_h \in W^{1,p}(\Om \setminus \Gamma_h)$ with $u_h \to u$ strongly in $L^1(\Om)$, for $h$ large enough we have
\begin{equation}
\label{mlm}
\int_\Om f(x,\nabla u)\,dx < \int_\Om f(x,\nabla u_h)\,dx +k\hs^1(\Gamma_h).
\end{equation}
In fact, it is not restrictive to assume that $(\nabla u_h)_{h \in \N}$ is bounded in $L^p(\Om,\R^2)$ and that $\hs^1(\Gamma_h) \le C$. By Ambrosio's lower semicontinuity theorem \cite{A2} we have 
$$
\int_\Om f(x,\nabla u)\,dx\le \liminf_{h \to +\infty} \int_\Om f(x,\nabla u_h)\,dx.
$$
If the sequence $(\hu(\Gamma_h))_{h \in \N}$ is not infinitesimal, then \eqref{mlm} clearly holds. If $\hu(\Gamma_h) \to 0$, we have $\hs^1(\Gamma_h)\le l^*$ for $h$ large enough, and hence \eqref{mlm} follows from \eqref{mmc}. 
\par
Actually, using Ambrosio's compactness theorem \cite{A2}, one can show that the local minimality in the $L^1$-topology and the minimality result of Theorem \ref{mainthm} are equivalent.
}
\end{remark}

In order to prove Theorem \ref{mainthm} we will use the main estimate given by Theorem \ref{mainestthm}.
%which is given by the analysis of the associated dual problem. 
Our aim is to construct a vector field $\tau \in L^q(\Om;\R^2)$ which is an admissible competitor in \eqref{mainest} and which shows that the difference between the bulk energies of $u$ and $u_\Gamma$ is smaller than $k\hu(\Gamma)$. In order to do so, we need some preliminary lemmas.

\begin{lemma}
\label{defeta}
Let $x \in \Om$, and let $r>0$ be such that $B_{2r}(x) \subseteq \Om$. 
Let $\varphi$ be a smooth function with $0 \le \varphi \le 1$, 
$\varphi=0$ on $B_r(x)$, $\varphi=1$ outside $B_{2r}(x)$ 
and $\|\nabla \varphi\|_{\infty} \le \frac{2}{r}$.
Then there exists $\eta \in L^q((B_{2r}(x)\setminus B_r(x));\R^2)$ with $q:=p'=p/(p-1)$ such that 
$$
\begin{cases}
\Div \, \eta=-\Div(\varphi \sigma)  & \text{on }B_{2r}(x) \setminus B_r(x), \\
\eta \cdot n=0  & \text{ on } 
\partial (B_{2r}(x)  \setminus B_{r}(x))
\end{cases}
$$ 
and
\begin{equation}\label{steta}
\int_{B_{2r}(x)\setminus B_r(x)}|\eta|^q\,dx \le C \int_{B_{2r}(x) \setminus B_r(x)}|\sigma|^q\,dx,
\end{equation}
where $\sigma:=\partial f(x,\nabla u)$ is the stress of the elastic solution $u$, $n$ is the outer normal to $\partial (B_{2r}(x)  \setminus  B_{r}(x))$, and $C$ is a constant independent of $r$.
\end{lemma}

\begin{proof}
Let us set  $\eta:=|\nabla v|^{p-2}\nabla v$, where $v \in W^{1,p}(B_{2r}(x)\setminus B_r(x))$  satisfies the equation
\begin{equation}
\label{defv}
\begin{cases}
\Div (|\nabla v|^{p-2}\nabla v)=-\Div(\varphi \sigma) &\text{in }B_{2r}(x)\setminus B_r(x),\\
|\nabla v|^{p-2}\nabla v \cdot n=0 &\text{on }\partial (B_{2r}(x)\setminus  B_r(x)).
\end{cases}
\end{equation}
Notice that the equation is well posed because
\begin{equation}
\label{divfisigma}
-\Div(\varphi \sigma)=-\nabla \varphi \cdot \sigma \in L^q(\Om)
\end{equation}
and
$$
\int_{B_{2r}(x) \setminus B_r(x)}\Div(\varphi \sigma)\,dx=
\int_{\partial B_{2r}(x)}\sigma \cdot n\,d\hu=0
$$
since $\sigma$ is divergence free in $\Om$.
\par
It remains to prove inequality \eqref{steta}. To this aim, note that we can always assume that $v$ has zero mean value, so that 
by Poincar\'e inequality and by a rescaling argument, we have that there exists $C>0$ independent on $r$ such that
\begin{equation}
\label{poincare2}
\int_{B_{2r}(x) \setminus B_r(x)}|v|^p\,dx \le
Cr^p
\int_{B_{2r}(x) \setminus B_r(x)}|\nabla v|^p\,dx. 
\end{equation}
Recalling that $\|\nabla \varphi\|_\infty \le 2/r $
and taking into account \eqref{divfisigma} and \eqref{poincare2} we get
\begin{multline*}
\int_{B_{2r}(x) \setminus B_r(x)}|\nabla v|^p\,dx=
\int_{B_{2r}(x) \setminus B_r(x)}\nabla \varphi \cdot \sigma v \,dx
\le \|\nabla \varphi \cdot \sigma\|_{L^q(B_{2r}(x) \setminus B_r(x))}
\|v\|_{L^p(B_{2r}(x) \setminus B_r(x))} \\
\le
\frac{1}{r} \|\sigma\|_{L^q(B_{2r}(x) \setminus B_r(x),\R^2)} C^{1/p}r\|\nabla v\|_{L^p(B_{2r}(x) \setminus B_r(x),\R^2)}\\=
C^{1/p}\|\sigma\|_{L^q(B_{2r}(x) \setminus B_r(x),\R^2)} \|\nabla v\|_{L^p(B_{2r}(x) \setminus B_r(x),\R^2)}
\end{multline*}
so that
$$
\int_{B_{2r}(x) \setminus B_r(x)}|\eta|^q\,dx = \int_{B_{2r}(x) \setminus B_r(x)}|\nabla v|^p\,dx \le C^{q/p} \int_{B_{2r}(x) \setminus B_r(x)}|\sigma|^q\,dx,
$$
and this concludes the proof.
\end{proof}

Now we need to construct a suitable field $\eta$ also around points $x$ on the boundary of $\Om$.  Since $\Om$ is Lipschitz, for every $x \in \partial \Om$ we can find an orthogonal coordinate system $(x_1',x_2')$ with origin at $x$, $\eps_1,\eps_2>0$ and a Lipschitz function $g:[-\eps_1,\eps_1] \to [-\eps_2,\eps_2]$ such that setting 
\begin{equation}
\label{defR}
R_r(x):=\{(x_1',x_2'): |x_1'| \le r\eps_1, |x_2'|\le r\eps_2\}
\end{equation}
we have for $r$ small enough
\begin{equation}\label{retsmal}
\Om \cap R_r(x)=\{(x_1',x_2') \in R_r(x)\,:\,x_2' \ge g(x_1')\}.
\end{equation}
Notice moreover that the Lipschitz constant $C_g$ of the function $g$ is determined only by $\Om$, and that we can assume $C_g \eps_1 < \eps_2$.
Let us set
\begin{equation}
\label{defAr}
A_r(x):=(R_{2r}(x) \setminus R_r(x)) \cap \Om.
\end{equation}

\begin{lemma}
\label{defetabdry}
Let $r>0$ be small enough (so that~\eqref{retsmal} holds everywhere on
$\partial\Om$). Let $x \in \partial \Om$ such that one of the
following three situations holds:
%% We distinguish three cases:
\begin{itemize}
\item[(1)] $x \in \partial_D \Om$ and $A_r(x)\cap\partial \Om\subset
\partial_D\Om$;
\item[(2)] $x \in \tint{\partial_N \Om}$ and
$A_r(x)\cap\partial \Om\subset \partial_N\Om$ (where $\tint{\cdot}$
  indicates the interior relative to $\partial \Om$);
\item[(3)] $x \in \partial \Om \setminus (\partial_D \Om \cup \tint{\partial_N \Om})$.
\end{itemize}
Let $\varphi$ be a smooth function with $0 \le \varphi \le 1$, $\varphi=0$ on $R_r(x)$, $\varphi=1$ outside $R_{2r}(x)$
and $\|\nabla \varphi\|_{\infty} \le \frac{2}{r}$.
Then there exists $\eta \in L^q(A_r(x);\R^2)$ with $q:=p'=p/(p-1)$ such that 
$$
\begin{cases}
\Div \, \eta=-\Div(\varphi \sigma)  & \text{in }A_r(x)\\
\eta \cdot n=0  & \text{on } \partial A_r(x) \cap (\Om \cup \partial_N \Om)
\end{cases}
$$
and
$$
\int_{A_r(x)}|\eta|^q\,dx \le C\int_{A_r(x)}|\sigma|^q\,dx,
$$
where $\sigma:=\partial f(x,\nabla u)$ is the stress of the elastic solution $u$, 
$n$ is the outer normal to $\partial A_r(x) \cap (\Om \cup \partial_N \Om)$,
and $C$ is a constant depending only on $\Om$.
\end{lemma}
We observe that only the points in $\partial\Om$ at distance less than
$r\sqrt{\eps_1^2+\eps_2^2}$ to $\partial \Om \setminus (\partial_D \Om
\cup \tint{\partial_N \Om})$ might not fall into one of first two cases.
\begin{proof}
In all the three cases of the lemma,
we will use the fact that the Poincar\'e inequality holds in
$W^{1,p}(A_r(x))$ with a constant that rescales as $r$, i.e., there
exists a positive constant $C$ independent of $r$ such that
\begin{equation}
\label{poincare}
\int_{A_r(x)}|v|^p\,dx \le
Cr^p
\int_{A_r(x)}|\nabla v|^p\,dx 
\end{equation}
for all $v \in W^{1,p}(A_r(x))$ with 
$$
\int_{A_r(x)}v=0
\qquad
\text{or}
\qquad
v=0 \text{ on }\partial_D \Om \cap \partial A_r(x).
$$
This can be seen rescaling $A_r(x)$ with the transformation $T_r(x',y')=(\frac{x'}{r},\frac{y'}{r})$,
and using Proposition~\ref{th:PK} (see Appendix~\ref{app:PK}) in the
domains $\{T_r(A_r(x))\}$:
it shows that
the Poincar\'e inequality holds in $T_r(A_r(x))$
with a constant that is independent of $r$, and by rescaling we deduce
that \eqref{poincare} holds.
\par
In case (1), we can consider $\eta:=|\nabla v|^{p-2}\nabla v$, with
$v \in W^{1,p}(A_r(x))$ satisfying the equation
\begin{equation}
\label{defv1}
\begin{cases}
\Div (|\nabla v|^{p-2}\nabla v)=-\Div(\varphi \sigma) &\text{in }A_r(x),\\
v=0 &\text{on }\partial A_r(x) \cap \partial_D \Om, \\
|\nabla v|^{p-2} \nabla v \cdot n=0 &\text{on } \partial A_r(x) \cap \Om.
\end{cases}
\end{equation}
Notice that $-\Div(\varphi \sigma)=-\nabla \varphi \cdot \sigma$ ($\sigma$ is divergence free) and $\|\nabla \varphi\|_\infty \le \frac{1}{r}$. Taking into account the Poincar\'e inequality \eqref{poincare} we get
\begin{multline*}
\int_{A_r(x)}|\nabla v|^p\,dx=
\int_{A_r(x)}\nabla \varphi \cdot \sigma v \,dx
\le \|\nabla \varphi \cdot \sigma\|_{L^q(A_r(x))}\|v\|_{L^p(A_r(x))} \\
\le
\frac{1}{r} \|\sigma\|_{L^q(A_r(x),\R^2)} C^{1/p}r\|\nabla v\|_{L^p(A_r(x),\R^2)}=
C^{1/p}\|\sigma\|_{L^q(A_r(x),\R^2)} \|\nabla v\|_{L^p(A_r(x),\R^2)}
\end{multline*}
so that
$$
\int_{A_r(x)}|\eta|^q\,dx=\int_{A_r(x)}|\nabla v|^p\,dx \le C^{q/p} \int_{A_r(x)}|\sigma|^q\,dx.
$$
\par
Cases (2) and (3) can be treated as case (1) considering $\eta:=|\nabla v|^{p-2}\nabla v$ with $v \in W^{1,p}(A_{r}(x))$ defined by the equations 
\begin{equation}
\label{defv2}
\begin{cases}
\Div (|\nabla v|^{p-2}\nabla v)=-\Div(\varphi \sigma) &\text{in }A_r(x),\\
|\nabla v|^{p-2} \nabla v \cdot n=0 &\text{on } \partial A_r(x),
\end{cases}
\end{equation}
and
\begin{equation}
\label{defv3}
\begin{cases}
\Div (|\nabla v|^{p-2}\nabla v)=-\Div(\varphi \sigma) &\text{in }A_r(x),\\
|\nabla v|^{p-2} \nabla v \cdot n=0 &\text{on } \partial A_r(x) \cap (\Om \cup \partial_N \Om), \\
v=0 &\text{on }\partial A_r(x) \cap \partial_D \Om
\end{cases}
\end{equation}
respectively. Notice in particular that equation \eqref{defv2} is well posed since its right hand side has zero mean value, because $\sigma$ is divergence free in $\Om$, and $\sigma \cdot n=0$ on $\partial_N \Om$ so that
$$
\int_{\Om}\Div(\varphi \sigma)\,dx=
\int_{A_r}\Div(\varphi \sigma)\,dx=
\int_{\partial A_r}\sigma \cdot n\,d\hu=0.
$$
\end{proof}

We are now in a position to prove our minimality result.

\begin{proof}[Proof of Theorem \ref{mainthm}]\label{proofmainthm}
We claim that there exist at most $m$ open balls $B_{r_1}(x_1),\dots,B_{r_k}(x_k)$, $k \le m$ such that $\Gamma \subseteq \cup_i B_{r_i}(x_i)$,
$$
r_i \le C\hu(B_{r_i}(x_i) \cap \Gamma),
$$
and 
$$
B_{2r_i}(x_i) \cap B_{2r_j}(x_j)=\emptyset
\qquad\text{for all }i \not=j,
$$
where $C$ depends only on $m$.
In fact let us consider the decomposition of $\Gamma$ in its connected components, i.e.,
$$
\Gamma:=\Gamma_1 \cup \dots \cup \Gamma_k,
$$
with $k \le m$. For all $i=1,\dots,k$ let $B_{s_i}(y_i)$ be an open ball with $s_i=\hu(\Gamma_i)$ and such that $\Gamma_i \subseteq B_{s_i}(y_i)$. 
If the balls $B_{2s_i}(y_i)$ are disjoint, then the covering $\{B_{s_i}(y_i)\}_{i=1,\dots,k}$ satisfies the claim. Otherwise we proceed in this way. Let us consider
$$
\bs^1:=\bigcup_{i=1}^{k} B_{2s_i}(y_i),
$$
and let $\bs^1_j$, $j=1,\dots,\tilde{k} \le k-1$ be its connected components. 
For all $j$, let $B_{\tau_j}(z_j)$ be an open ball with $\tau_j=\diam{\bs^1_j}$ such that 
$\bs^1_j \subseteq B_{\tau_j}(z_j)$. Again, if the balls $B_{2\tau_j}(z_j)$ are disjoint, then the covering $\{B_{\tau_j}(z_j)\}_{j=1,\dots,\tilde k}$ satisfies the claim.
Otherwise we construct in a similar way as before the set $\bs^2$ which has at most $k-2$ connected components. Clearly in at most $m$ steps we come up with at most $m$ balls satisfying the requirements of the claim.
\par
Since $\Om$ is Lipschitz, we have that if $\hu(\Gamma)$ is
sufficiently small (depending on $\Om$), we can assume that the balls
$\{B_{r_i}(x_i)\}_{i=1,\dots,k}$ intersecting $\partial \Om$ can be
replaced by rectangles $R_{s_j}(y_j)$ of the form \eqref{defR} satisfying
\eqref{retsmal}, centered at some point $y_j$ that falls into case (1), (2)
or (3) of Lemma~\ref{defetabdry}. More precisely, 
thare exist a constant $C$ depending only on $\Om$ and $m$, and
there exist at most $m$ open balls $\{B_{r_i}(x_i)\}$ and at most $m$ rectangles $\{R_{s_j}(y_j)\}$ defined in \eqref{defR} with $y_j$ falling into case (1), (2) or (3) of Lemma~\ref{defetabdry}, such that $B_{2r_i}(x_i) \subseteq \Om$, and
$$
\Gamma \subseteq \cup_{i,j} B_{r_i}(x_i) \cup R_{s_j}(y_j),
$$
\begin{equation}
\label{estdiam}
\diam{B_{r_i}(x_i)} \le C \hu(B_{r_i}(x_i) \cap \Gamma),
\qquad
\diam{R_{s_j}(y_j)} \le C\hu(R_{s_j}(y_j) \cap \Gamma),
\end{equation}
$$
B_{2r_i}(x_i) \cap B_{2r_j}(x_j)=\emptyset,
\qquad
R_{2s_i}(y_i) \cap R_{2s_j}(y_j)=\emptyset,
\qquad
B_{2r_i}(x_i) \cap R_{2s_j}(y_j)=\emptyset
\qquad\text{for all }i,j.
$$
\par
Let $\varphi$ be a smooth function with 
$$
0 \le \varphi \le 1,
\qquad 
\varphi=0 \text{ on }\bigcup_{i,j} \left( B_{r_i}(x_i) \cup R_{s_j}(y_j) \right),
\qquad
\varphi=1 \text{ outside }\bigcup_{i,j} \left( B_{2r_i}(x_i) \cup R_{2s_j}(y_j) \right).
$$
Let us denote with $A_{r_i}(x_i)$ and $A_{s_j}(y_j)$ the sets $B_{2r_i}(x_i) \setminus B_{r_i}(x_i)$ and $R_{2s_j}(y_j) \setminus R_{s_j}(y_j)$ respectively.
Let 
$$
\eta_i \in L^q\left(A_{r_i}(x_i);\R^2\right)
\quad\text{and}\quad
\eta_j \in L^q\left(A_{s_j}(y_j);\R^2\right)
$$ 
with $q:=p'=p/(p-1)$ be the vector fields given by Lemmas \ref{defeta} and \ref{defetabdry}.
Let us consider $\eta \in L^q(\Om;\R^2)$ defined as
$$
\eta:=
\begin{cases}
\eta_i      & \text{in }A_{r_i}(x_i), \\
\eta_j      & \text{in }A_{s_j}(y_j), \\
0             & \text{otherwise,}
\end{cases}
$$
and let us set
$$
\tau:= \varphi \sigma+ \eta,
$$
where $\sigma=\partial f(x,\nabla u)$ is the stress of the elastic solution $u$. By construction we have that 
\begin{equation}
\label{tauok}
\int_\Om \tau \cdot \nabla v=0 \qquad \text{for all }v \in \as(\Gamma),
\end{equation}
where $\as(\Gamma)$ is as in \eqref{agamma}. Moreover we have that for all $i,j$
$$
\int_{A_h} |\eta|^q\,dx \le 
C\int_{A_h} |\nabla u|^p\,dx,
$$
where $A_h$ denotes one of the $A_{r_i}(x_i)$'s or one of the $A_{s_j}(y_j)$'s,
and $C$ is a constant depending only on $\Om$ and the bulk energy density $f$.
\par
In view of \eqref{tauok} we can put $\tau$ in inequality \eqref{mainest} getting
\begin{equation}
\label{mainesttaun}
\int_\Om [f(x,\nabla u)-f(x,\nabla u_{\Gamma})]\,dx \le \sum_j 
\int_{A_h}
[\tau -\sigma] \cdot [\partial f^*(x, \tau)-\partial f^*(x,\sigma)]\,dx.
\end{equation}
Since $f^*$ is the convex conjugate of $f$, and since $f$ satisfies the growth conditions
\eqref{fcrescita}, we deduce that
$$
|f^*(x,\zeta)| \le C(|\zeta|^q+1)
$$
and
$$
|\partial_\zeta f^*(x,\zeta)| \le C(|\zeta|^{q-1}+1).
$$
We claim that for every $A_h$
\begin{equation}
\label{ineqdiam}
\int_{A_h}
[\tau -\sigma] \cdot [\partial f^*(x, \tau)-\partial f^*(x,\sigma)] \le C \diam{A_h}^\alpha,
\end{equation}
where $C$ is independent of $A_h$, and depends only on $m$, $\Om$ and $f$. Here $\alpha>1$ is the exponent defining the weak singularities of $u$ (see Definition \ref{weaksingdef}). From \eqref{ineqdiam} and \eqref{mainesttaun}, taking into account \eqref{estdiam} and the fact that $\alpha>1$, we obtain that there exists $l^*$ depending only on $\Om$, $m$, $f$ and $k$ such that for every $\Gamma \in \ks_m(\Omb)$ with $\hu(\Gamma)<l^*$ we have
$$
\int_\Om [f(x,\nabla u)-f(x,\nabla u_{\Gamma})]\,dx \le C \sum_h \diam{A_h}^\alpha\le
C\hu(\Gamma)^\alpha<
k\hu(\Gamma),
$$
so that the minimality result holds.
\par
In order to conclude the proof, we need to show that claim \eqref{ineqdiam} holds true.
This can be seen making all the products and estimating each addend. Let us check the first one, the other ones being similar. We have
$$
\int_{A_h} \tau \cdot \partial f^*(x, \tau) \,dx \le
C\int_{A_h} (|\sigma|+|\eta|)\cdot(|\sigma|^{q-1}+|\eta|^{q-1}+1)\,dx.
$$
Then, in view of Lemmas \ref{defeta} and \ref{defetabdry}, since $u$ has uniformly weak singularities in $\Om$ we get for $r$ small enough
$$
\int_{A_h} |\sigma|^q \,dx \le
C\int_{A_h} |\nabla u|^p\,dx \le C \diam{A_h}^\alpha,
$$
$$
\int_{A_h} |\eta|^q \,dx \le
C\int_{A_h} |\nabla u|^p\,dx \le C \diam{A_h}^\alpha,
$$
$$
\int_{A_h} |\eta| |\sigma|^{q-1} \,dx \le
\left( \int_{A_h} |\eta|^q\,dx \right)^{\frac{1}{q}} 
\left( \int_{A_h} |\sigma| ^q\,dx \right)^{\frac{1}{p}}
\le  C \diam{A_h}^\alpha,
$$
$$
\int_{A_h} |\eta|^{q-1} |\sigma| \,dx \le
\left( \int_{A_h} |\eta|^q\,dx \right)^{\frac{1}{p}} 
\left( \int_{A_h} |\sigma|^q\,dx \right)^{\frac{1}{q}}
\le  C \diam{A_h}^\alpha,
$$
$$
\int_{A_h} |\sigma| \,dx \le
C\left( \int_{A_h} |\sigma|^q\,dx \right)^{\frac{1}{q}}
\diam{A_h}^{\frac{2}{p}}
\le C \diam{A_h}^\alpha,
$$
and
$$
\int_{A_h} |\eta| \,dx \le
C\left( \int_{A_h} |\eta|^q\,dx \right)^{\frac{1}{q}}
\diam{A_h}^{\frac{2}{p}}
\le C \diam{A_h}^\alpha.
$$
Summing up we obtain that \eqref{ineqdiam} holds.
\end{proof}

\begin{remark}
\label{hypsharp}
{\rm
Notice that the arguments of the previous proof also work in the case in which the elastic solution $u$ has critical singularities, i.e. in the case $\alpha=1$, provided that
there exists $\eps>0$ and $\delta$ sufficiently small (depending only on $\Om$, $m$, $f$ and $k$), such that
$$
\int_{B_r(x) \cap \Omb}|\nabla u|^p\,dx\le \delta r
\qquad
\text{for every }r<\eps.
$$
}
\end{remark}

Using the same arguments of the proof of Theorem \ref{mainthm}, we can easily deduce
the following localized version of the minimality result.

\begin{proposition}
\label{mainthmrmk}
Let the elastic solution $u$ of problem \eqref{elasticsol} have only uniformly weak singularities in $A$, where $A$ is an open subset of $\Om$. Then there exists a critical length $l^*>0$ depending only on $A$, $m$, $f$, and $k$ such that for all $\Gamma \in \ks_m(\Omb)$ with $\Gamma \subseteq \bar A$ and $\hu(\Gamma)<l^*$ we have
$$
\int_\Om f(x,\nabla u)\,dx<\int_\Om f(x,\nabla u_\Gamma)\,dx+k\hu(\Gamma),
$$
where $u_{\Gamma}$ is a minimum of \eqref{elasticsol}.
\end{proposition}

\begin{remark}
\label{strongfalse}
{\rm {\bf (The case of strong singularities)}
Theorem \ref{mainthm} is false if the elastic solution $u$ has {\it strong singularities} in $\Om$, i.e., there exists $x \in \Omb$ such that
\begin{equation}
\label{strongsing}
\limsup_{r \to 0} \frac{1}{r}\int_{B_r(x) \cap \Om}|\nabla u|^p\,dx=+\infty.
\end{equation}
In fact if such a point $x$ exists, then the pair $(v_r,\partial B_r(x) \cap \Omb)$, where $v_r$ is defined as
\begin{equation}
\label{vr}
v_r(y):=
\begin{cases}
u(y)      & \text{in }\Omb \setminus B_r(x), \\
0      & \text{in }\Omb \cap B_r(x)
\end{cases}
\end{equation}
is energetically more convenient with respect to $u$ for $r$ small.
Note also that this example needs the right hand-side of
\eqref{strongsing} to be just greater than $\frac{2\pi}{k}$}.
\end{remark}

The previous remark and Proposition \ref{mainthmrmk} seem to suggest that in presence of strong singularities, energetically convenient small cracks prefer to stay near the singular points. Let us prove that this intuition is indeed true in the case in which $u$ has only one point of strong singularity $x$. 
\par
For every $l>0$ let $\Gamma_l \in \ks_m(\Omb)$ be such that $(u_{\Gamma_l},\Gamma_l)$ minimizes the total energy \eqref{totalenergy} among all pairs $(u_\Gamma,\Gamma)$ with $\Gamma \in \ks_m(\Omb)$ and $\hu(\Gamma) \le l$. The existence of $\Gamma_l$ can be proved using the direct method of the Calculus of Variations in view of the lower semicontinuity of the $\hu$-measure with respect to Hausdorff converging sequences in $\ks_m(\Omb)$ given by G\c{o}\l ab's Theorem (see for example \cite{DMT} for details). 
\par
Notice that $\Gamma_{l}\neq \emptyset$ for $l$ small because
the pair $(v_l, \partial B_l)$, where $v_l$ is defined in \eqref{vr}, is energetically more convenient than the elastic solution. Moreover, for every $r>0$, we have that
$\Gamma_l \cap B_r(x) \not= \emptyset$ when $l$ is small enough
because otherwise, in view of Proposition \ref{mainthmrmk}, $(u_{\Gamma_l},\Gamma_l)$ would not be energetically more convenient than the elastic solution. So we deduce that the following proposition holds.
 
\begin{proposition}
\label{smcr}
Assume that $x \in \Om$ is a point of strong singularity for the elastic solution $u$, and that $u$ has uniformly weak singularities in $\Omb \setminus B_r(x)$ for every $r>0$. Then for every neighborhood $U$ of $x$, if $l$ is small enough
we have $\Gamma_{l} \neq \emptyset$, and $\Gamma_l \cap U \not= \emptyset$.
\end{proposition}

\begin{remark}
\label{overviewsec}
{\rm {\bf (Singularities in materials)}
Theorem \ref{mainthm} and Remark \ref{strongfalse} show that the quantity
\begin{equation}
\label{gradball}
\int_{\Om \cap B_r(x)}|\nabla u|^p\,dx,
\end{equation}
where $u$ is a solution of
\begin{equation}
\label{eqelastic}
\begin{cases}
\Div \partial_\xi f(x,\nabla u)=0 &\text{in }\Om \\
u=\psi &\text{on }\partial_D \Om,
\end{cases}
\end{equation}
determines the feasibility of the appearance of small cracks in $\Omb$ when imposing a boundary datum $\psi$.
\par
Let us recall some known result concerning the behavior of \eqref{gradball}.
As before, we assume that $f: \Om \times \R^2 \to \R$ is a Carath\'eodory function satisfying \eqref{fstruttura} and  such that for a.e. $x \in \Om$ and for all $\xi \in \R^2$
$$
\lambda |\xi|^p \le f(x,\xi) \le \Lambda |\xi|^p,
$$
where $\lambda, \Lambda>0$. 
\par
By Theorem \ref{mainthm}, small cracks are energetically convenient if $u$ has uniformly weak singularities in $\Omb$ according to Definition \ref{weaksingdef}. This is certainly the case if $\nabla u$ has a summability sufficiently higher than $p$, as one can check by means of H\"older inequality. In fact we have
\begin{equation*}
\int_{B_r(x)} |\nabla u|^{p}\,dx
\le \left( \int_{B_r(x)} |\nabla u|^q \,dx \right)^{p/q} |B_r(x)|^{(q-p)/q}
\le \left( \int_{\Om} |\nabla u|^q \,dx \right)^{p/q} r^{2(q-p)/q}.
\end{equation*}
We deduce that $u$ has weak singularities in $\Omb$ if $\nabla u \in L^q(\Om;\R^2)$ with $q>2p$.
\par
It is well known that if $\Om$, $\psi$ and $f$ are sufficiently regular, then the solution $u$ of \eqref{eqelastic} is regular, so that $u$ has only uniformly weak singularities in $\Om$. 
However the assumption of regularity of $f$ with respect to the variable $x$ is not suitable for applications to continuum mechanics, since discontinuity in the variable $x$ models the important case of composite materials. 
\par
Several papers in the literature address the issue of higher integrability properties of the gradient of solutions of \eqref{eqelastic} without continuity assumptions on $x$ (see \cite{CP}, \cite{LiNir}, \cite{LiVo}, \cite{LN} and references therein), namely properties of type $\nabla u \in L_{\rm loc}^q(\Om;\R^2)$ for some $q>p$. The restriction of the behavior of $\nabla u$ to compactly contained open subsets of $\Om$ is the natural price to pay in order to concentrate on properties depending only on the material (i.e., on the bulk energy density $f$) and not on the boundary datum. 
\par
For example a regularity result by Caffarelli and Peral \cite[Theorem C]{CP} implies that, under mild assumptions on $f$, if for all $q>p$ there exists $\eps=\eps(q)>0$ such that if
$$
||\xi|^{p-2}\xi-\partial_\xi f(x,\xi)| \le \eps |\xi|^{p-1}
$$
(i.e., $\partial_\xi f$ is sufficiently near to the $p$-Laplacian operator), then $u \in W^{1,q}_{\rm loc}(\Om)$. In particular
if $\eps$ is small enough, $u$ has uniformly weak singularities on every compactly contained open subset of $\Om$. 
\par
In the linear case $p=2$, the hypothesis of the regularity result by Caffarelli and Peral \cite{CP} reduces to the requirement that
$$
\Lambda/\lambda<1+\eps
$$
for $\eps$ sufficiently small. However Leonetti and Nesi \cite{LN} found the optimal integrability exponent for $\nabla u$ in terms of  $K:= (\Lambda/\lambda)^{1/2}$, namely they proved that 
$$
\nabla u \in L^s_{\rm loc}(\Om)
\quad\text{for every }s<2K/(K-1).
$$
So in order to guarantee that $u$ has uniformly weak singularities in compactly contained open subsets of $\Om$, it suffices that $\Lambda/\lambda<4$. 
\par
A famous example due to Meyers \cite{M} shows that if $\Lambda/\lambda=4$
strong singularities may appear inside the body. In fact, assuming that the origin belongs to $\Om$
we can consider
$$
A(x):= 1/K  n \otimes n + K \tau \otimes \tau,
$$
where $n:=\frac{x}{|x|}$, $\tau$ is obtained from $n$ through a rotation of 90 degrees counterclockwise, and $a \otimes b$ denotes the matrix with coefficients $(a \otimes b)_{ij}=a_ib_j$.  Then it is easy to see that the solution of  
$$
\begin{cases}
\Div(A(x)\nabla u)=0 &\text{in }\Om \\
u(x)=x_1 &\text{on }\partial \Om
\end{cases}
$$
is given by
$$
u(x):=|x|^{K-1} x_1.
$$
A simple computation shows that if $K\ge 2$, i.e., $\Lambda/\lambda \ge 4$, the origin is a point of strong singularity for $u$.
}
\end{remark}

\begin{remark}
\label{secvectorial}
{\rm {\bf (The vectorial $2D$-case)}
The minimality result given by Theorem \ref{mainthm} holds also in the case of 
$\R^M$-valued displacements provided that we choose them in 
the Deny-Lions space
\begin{equation}
\label{denylions}
L^{1,p}(\Om \setminus \Gamma;\R^M):=\{v \in W^{1,p}_{\rm loc}(\Om \setminus \Gamma,\R^M)\,:\, \nabla v \in L^p(\Om \setminus \Gamma; M^{M \times 2})\}.
\end{equation}
We have that $W^{1,p}(\Om \setminus \Gamma,\R^M) \subseteq L^{1,p}(\Om \setminus \Gamma;\R^M)$, and the two spaces are equal if $\Om \setminus \Gamma$ is sufficiently regular (for example an union of a finite number of Lipschitz domains). A notion of trace for functions in $L^{1,p}(\Om \setminus \Gamma; \R^M)$ near the points of $\partial_D \Om \setminus \Gamma$ is well defined, so that given $\psi \in W^{1,p}(\Om,\R^M)$ and $\Gamma \in \ks_m(\Omb)$, the displacement $u_\Gamma$ is a solution of the minimum problem
\begin{equation}
\label{elasticsolvectorial}
\min \left\{ \int_\Om f(x,\nabla u)\,dx\,:\, u \in L^{1,p}(\Om \setminus \Gamma;\R^M), 
u=\psi \text{ on }\partial_D \Om \setminus \Gamma \right\}.
\end{equation}
(Notice that in this vectorial setting the maximum principle does not hold, so that only a control on the gradient is available and this is why Deny-Lions spaces are required).
Since Definition \ref{weaksingdef} relies only on the behavior of $\nabla u$, it turns out that the notion of uniformly weak singularities is well defined in the context of $L^{1,p}$-spaces.
\par
The minimality result in the vectorial setting follows because estimate \eqref{mainest} still holds provided we set
$$
\as(\Gamma):= \{v \in L^{1,p}(\Om \setminus \Gamma;\R^M)\,:\, v=0 \text{ on }\partial_D \Om \setminus \Gamma \},
$$
and the constructions of Lemmas \ref{defeta} and \ref{defetabdry} can easily be adapted to the case of matrix valued vector fields. 
}
\end{remark}

\begin{remark}
\label{higherdim}
{\rm {\bf (The $N$-dimensional case)}
Let us consider the case $\Om \subseteq \R^N$ with $N \ge 3$. The two dimensional setting is employed in the proof of Theorem \ref{mainthm} only to ensure the existence of a covering of the crack $\Gamma$ which satisfies conditions \eqref{estdiam}. In the case $\Gamma$ is connected, the covering condition can be reduced to the existence of an open set $A$ (a ball if $\Gamma$ is well inside $\Om$, or a rectangle if it is near the boundary) such that $\Gamma \subseteq A$ and
$$
\diam{A} \le C \hu(\Gamma),
$$
where $C$ is a given constant. This inequality is not implied by connectedness in dimension $N \ge 3$, because ``{\it needle-like}'' cracks can have very small $\hn$-measure and big diameter. So the machinery of the proof of Theorem \ref{mainthm} can be employed in dimension $N \ge 3$ if we restrict to a class of admissible cracks which excludes the elongated  ones. Namely we can consider the family
\begin{equation}
\label{Ncracks}
\ks^C(\Omb):=\{\Gamma \subseteq \Omb\,:\, \Gamma \text{ is closed and } \diam{\Gamma} \le C\hn(\Gamma)\},
\end{equation}
where $C$ is a given constant. The displacement $u_\Gamma$ relative to 
a crack $\Gamma \in \ks^C(\Omb)$ and the boundary datum $\psi \in W^{1,p}(\Om,\R^M)$, with $M \ge 1$, is given again by problem 
\begin{equation}
\label{elasticsolNdim}
\min \left\{ \int_\Om f(x,\nabla u)\,dx\,:\, u \in L^{1,p}(\Om \setminus \Gamma;\R^M), 
u=\psi \text{ on }\partial_D \Om \setminus \Gamma \right\},
\end{equation}
where the space $L^{1,p}(\Om \setminus \Gamma;\R^M)$ is defined in \eqref{denylions}.
\par
The notion of uniformly weak singularities can be rephrased in the $N$-dimensional setting in the following way: we say that $u \in L^{1,p}(\Om; \R^N)$ has uniformly weak singularities in $\Omb$ if there exist constants $N-1<\alpha\le N$ and $C>0$ such that for every $x \in \Omb$ and for $r$ small enough
\begin{equation}
\label{weaksingNdim}
\int_{B_r(x) \cap \Om}|\nabla u|^p\,dx \le Cr^\alpha.
\end{equation}
The intuitive meaning of condition \eqref{weaksingNdim} is the same as in the two-dimensional setting, namely that the bulk energy of $u$ inside a ball of center $x$ and radius $r$ is asymptotically negligible for $r \to 0$ and uniformly in $x$ with respect to the surface of the ball.
}
\end{remark}

\section{The minimality result in linearized elasticity}
\label{planarelasticitysec}

We briefly show in this section that the results obtained up to now are
also true if the bulk energy $f(x,\nabla u)$ is replaced with a quadratic
linearized elasticity energy
$$
C(x)e(u):e(u)=C_{i,j,k,l}(x)e(u)_{i,j}e(u)_{k,l},
$$
that is a positive-definite quadratic form of the symmetrized gradient
$$
e(u):=(\nabla u + (\nabla u)^T)/2
$$ 
of the displacement $u:\Om\to \R^N$.
As before, the main conclusions will be drawn in dimension 2
(Theorem~\ref{mainthm} or rather, since $u$ is vectorial,
the result of Remark \ref{secvectorial}), whereas in higher dimension the same
restrictions of Remark \ref{higherdim} will apply.
\par
Let us consider in $\Om$ a measurable function $C$, such that for any
$x\in \Om$, $C(x)$ is a  $N\times N\times N\times N$ 4th-order tensor
that defines a positive-definite quadratic form on the vector
space of symmetric $N\times N$ matrices, that we denote by
$\SymN$. We assume that
for any $\xi\in \SymN$ and a.e. $x\in \Om$, it holds
\begin{equation}\label{elascoer}
\lambda |\xi|^2 \,\leq \, C(x)\xi:\xi \,\leq\, \Lambda|\xi|^2\,,
\end{equation}
where $\xi:\eta = \Tr(\xi \eta^T) = \xi_{i,j}\eta_{i,j}$ and $|\xi|^2=\xi:\xi$ is the standard Euclidean (Frobenius) norm.
\par
Given $\Gamma\subset \Omb$ a compact one-dimensional fracture,
the space of admissible displacements with finite energy will be
the space of measurable displacements $u:\Om\to \R^N$ whose symmetrized
distributional gradient in $\Om\setminus \Gamma$, denoted by $e(u)$,
is in $L^2(\Om\setminus \Gamma;\SymN)$, and that satisfy in some sense
$u=\psi$ on $\partial_D\Om\setminus\Gamma$. Thanks to 
\emph{Korn's inequality}, it is known that such a displacement
belongs in fact to $H^1_{\rm loc}(\Om\setminus \Gamma)$, and since
we have assumed that the boundary of $\Om$ is Lipschitz, we also have
that $u\in H^1(\Om\cap B)$  for any ball $B$ with $\overline{B}\cap
\Gamma=\emptyset$, so that the trace of $u$ on $\Om\setminus \Gamma$
is well defined. As in Section~\ref{secvectorial}, we introduce
the space
\[
LD^{1,2}(\Om\setminus \Gamma):=
\{ u\in H^1_{\rm loc} (\Om\setminus\Gamma;\R^N)\,:\, e(u)\in L^2(\Om\setminus
\Gamma;\SymN)\}
\]
and for $\Gamma$ a closed subset of $\Omb$,
the displacement $u_\Gamma$ is given by
\[
\min\left\{ \int_{\Om} C(x)e(u):e(u)\,dx\,:\, u\in LD^{1,2}(\Om\setminus \Gamma), u=\psi \textrm{ on } \partial_D\Om\setminus \Gamma\right\}\,.
\]
We denote by $u$ the solution  in the non-cracked domain.
The set of admissible variations is now
\[
\as(\Gamma):= \{v \in LD^{1,2}(\Om \setminus \Gamma)\,:\, v=0 \text{ on }\partial_D \Om \setminus \Gamma \}.
\]
The proof of Theorem~\ref{mainestthm} can be reproduced in this
situation, yielding the estimate (we assume $|\Gamma|=0$)
\[
\int_\Om [C(x)e(u):e(u)-C(x)e(u_\Gamma):e(u_\Gamma)]\,dx
\,\leq\,
2 \int_\Om [\tau-\sigma]:[C(x)^{-1}(\tau-\sigma)]\,dx,
\]
with $\sigma(x):=C(x)e(u)(x)$ and $\tau\in L^2(\Om;\SymN)$
is any stress field compatible
with the variations in $\as(\Gamma)$, that is, such that
\begin{equation}
\label{divgammasym}
\int_{\Om\setminus \Gamma} \tau:e(v)\,dx \,=\,0
\end{equation}
for any $v\in \as(\Gamma)$.
The last estimate, together with~\eqref{elascoer}, yields 
\begin{equation}
\label{estimelas}
\int_\Om [C(x)e(u):e(u)-C(x)e(u_\Gamma):e(u_\Gamma)]\,dx
\,\leq\, \frac{2}{\lambda} \int_\Om |\tau -\sigma|^2\,dx\,,
\end{equation}
for any $\tau$ satisfying \eqref{divgammasym}.
\par
In order to prove Theorem \ref{mainthm} in this setting, we do exactly
the same construction as before. We build a $\tau$ from $\sigma$
by letting $\tau=\sigma$ in $\Om$ except in a finite union of balls
or rectangles $B_{2r_i}(x_i)$ or $R_{2s_j}(y_j)$
(cf the proof of Theorem~\ref{mainthm} pp.~\pageref{proofmainthm} and following).
Inside the smaller balls/rectangles $B_{r_i}(x_i)$ and $R_{s_j}(y_j)$,
we choose $\tau=0$, and in each crown $A_{r_i}(x_i)=B_{2r_i}(x_i)\setminus B_{r_i}(x_i)$ or $A_{s_j}(y_j)=R_{2s_j}(y_j)\setminus R_{s_j}(y_j)$, $\tau$ is of the form: $\tau=
\varphi \sigma+\eta$ where $\varphi$ is the appropriate cut-off function.
\par
Again, to achieve \eqref{divgammasym}, one needs to choose
$\eta$ in an appropriate way. An additional difficulty here follows
from the fact that $\eta(x)$ has to be almost everywhere a $N\times N$
symmetric matrix. In order to find a suitable
$\eta$, one needs to replace problems~\eqref{defv},
\eqref{defv1}, \eqref{defv2} or~\eqref{defv3} by the appropriate variant.
\par
The right way to do is obviously to solve the vectorial equation
\begin{equation}\label{defvlinest}
\Div\, e(v)\ =\ -\Div\, (\varphi\sigma)
\end{equation}
in the appropriate domain, and to replace the Neumann boundary condition,
when present, with the corresponding condition $e(v)\cdot n=0$. We then
set $\eta=e(v)$ in each crown.
\par
We get an estimate on $\int_\Om |\tau -\sigma|^2\,dx$ from standard
estimates on $e(v)$, that will follow from  appropriate 
\textit{Poincar\'e-Korn} inequalities. For instance, one shows 
(see Appendix~\ref{app:PK} for details) that
\begin{equation}\label{poincarekorn}
\int_{A_r(x)} |v|^2\,dx \ \leq\ Cr^2\int_{A_r(x)} |e(v)|^2\,dx
\end{equation}
for any $v\in H^1(A_r(x);\R^N)$ with
\[
\left(\int_{A_r(x)} v(y)\,dy=0\ \textrm{ and }\ 
\int_{A_r(x)} y\times v(y)\,dy=0\,\right), 
\textrm{ or }\ 
v=0\textrm{ on } \partial_D\Om\cap\partial A_r(x)\,.
\]
The first set of conditions ensure that the ``average rigid motion''
of $v$ vanishes in $A_r(x)$. Rigid motions, of the form $x\mapsto
Kx+p$ with $K$ antisymmetric, are the kernel of the symmetrized
gradient (in any connected domain).
\par
Multiplying~\eqref{defvlinest} by $v$ and integrating by parts, we find
that
\begin{multline*}
\int_{A_r(x)} |e(v)|^2\,dx
 \ =\ -\int_{A_r(x)} \sigma:(\nabla\varphi\otimes v(x))\,dx
 \\
 \leq\ 
 \|\sigma\|_{L^2(A_r(x))} \|\nabla\varphi\|_{L^\infty(A_r(x))}
\left( \sqrt{C}r \|e(v)\|_{L^2(A_r(x))}\right)\,,
\end{multline*}
where we have used~\eqref{poincarekorn}. Since $\|\nabla \varphi\|_{L^\infty(A_r(x)}\le 2/r$ and $\eta=e(v)$, we deduce that
\[
 \int_{A_r(x)} |\eta(x)|^2\,dx\ \leq\ C\int_{A_r(x)}|\sigma(x)|^2\,dx
\]
for some constant $C$ that does not depend on $x$ or $r$. We conclude as in
the proof of Theorem~\ref{mainthm}.

\section{Qualitative properties of crack initiation}
\label{inisec}
In this Section we use the results of Section~\ref{mainresultsec} to address the problem of crack initiation in elastic bodies. We restrict our analysis to the case of antiplane elasticity. In view of the minimality result of Section~\ref{planarelasticitysec}, the same conclusions hold also for the case of planar linearized elasticity.
\par
First of all we consider the classical Griffith's theory
of quasistatic crack propagation, and we prove that it cannot explain the formation of a crack in an elastic body $\Om$ without singularities within the class $\ks_m(\Omb)$ of closed sets with a finite number of connected components, which is much richer than the family of smooth curves usually considered in the literature. 
\par
In the second part of the section, we remove the assumption, implicit in Griffith's theory, that the crack growth is progressive, namely that the length of the crack is continuous in time. Inspired by the variational theory of quasistatic crack evolution proposed by Francfort and Marigo in \cite{FM}, we replace the classical Griffith's equilibrium condition with a {\it static equilibrium condition} and an {\it energy balance}. The static equilibrium condition is a {\it unilateral minimality property} which states that, during the crack evolution, the total energy is minimal among all configurations with larger cracks (so that discontinuities of the crack's length are allowed). The energy balance requires that the total energy of the system evolves in relation with the power of external loads in such a way that no dissipation occurs. 
\par
Within this framework, we prove that a crack appears immediately at a point of strong singularity for the body. Moreover we prove that if the body has uniformly weak singularities, then it deforms elastically until a critical time $t_i$ after which a ''big'' crack $\Gamma(t)$ appears. These results have been established by Francfort and Marigo in \cite[Proposition 4.19, point (ii)]{FM} under the assumptions that the crack $\Gamma(t)$ is union of $m$ fixed curves $\{\gamma_i(t)\}_{i=1,\dots,m}$ which can be parametrized by arc length.
Thank to our local minimality result (Theorem \ref{mainthm}), we prove these facts removing the restrictions on the path of the crack. 
\par
The mathematical setting we consider is that of Section~\ref{mainresultsec}.
Namely $\Om$ is a bounded Lipschitz open set in $\R^2$, $\partial_D \Om \subseteq \partial \Om$ is open in the relative topology, and $\partial_N \Om:=\partial \Om \setminus \partial_D \Om$ is composed of a finite number of connected components.
The family of admissible cracks is given by the class $\ks_m(\Omb)$ 
of closed subsets of $\Omb$ with at most $m$ connected components and with finite length, while the class of admissible displacements relative to a crack $\Gamma$ is given by $W^{1,p}(\Om \setminus \Gamma)$ with $p \in ]1,+\infty[$. The total energy is given by
\begin{equation}
\label{griffithenergy}
\Es(u,\Gamma):= \int_\Om f(x,\nabla u)\,dx+k \hu(\Gamma),
\end{equation}
where $f$ is a Carath\'eodory function satisfying \eqref{fstruttura}, \eqref{fzero} and \eqref{fcrescita}, and $k>0$. 

\subsection{Crack initiation and Griffith's theory}
Let $\Gamma_0$ be a crack inside $\Omb$ of length $l_0$, and suppose that a boundary displacement $\psi$ is assigned on $\partial_D \Om \setminus \Gamma$. According to Griffith theory, $\Gamma_0$ is in equilibrium if, taken any family of increasing cracks $\Gamma_l$ containing $\Gamma_0$ with length $l_0+l$, then
\begin{equation}
\label{griffitheq}
\limsup_{l \to 0^+} \frac{\ws(l_0)-\ws(l_0+l)}{l} \le k,
\end{equation}
where $\ws(l_0)$ and $\ws(l_0+l)$ denote the bulk energy of the displacements $u_{\Gamma_0}$ and $u_{\Gamma_l}$ associated to the boundary datum $\psi$ and the cracks $\Gamma_0$ and $\Gamma_l$ respectively, and $k$ represents the toughness of the material. Moreover, during a {\it quasistatic crack evolution}, if $\Gamma_0$ propagates along the $\Gamma_l$, then \eqref{griffitheq} holds with equality.
\par
Let us prove that the rate of energy release that appear in the left hand side of \eqref{griffitheq} is zero in the case in which $\Gamma_0=\emptyset$ and the elastic solution $u$ relative to the boundary displacement $\psi$ has uniformly weak singularities. This means that the elastic configuration is always in equilibrium according to Griffith's theory, and moreover that a quasistatic crack evolution which begins in the elastic configuration remains at all subsequent times in the elastic regime, i.e., Griffith's theory cannot explain crack initiation.
\par
To this aim for every $l>0$, let us set
$$
\ws(l):=\inf \left\{ \int_\Om f(x,\nabla u_\Gamma)\,dx\,:\, \Gamma \in \ks_m(\Om),
\hu(\Gamma) \le l \right\},
$$
where $u_\Gamma \in W^{1,p}(\Om \setminus \Gamma)$ denotes the displacement associated to $\Gamma$ and the boundary datum $\psi$. Notice that we have clearly $\ws(0)=\int_\Om f(x,\nabla u)\,dx$, where $u$ is the elastic displacement relative to 
$\psi$.
\par
The following proposition holds.

\begin{proposition}
\label{griffith}
Let us assume that the hypothesis of Theorem \ref{mainthm} are fulfilled. Then we have
$$
\lim_{l \to 0^+}\frac{\ws(0)-\ws(l)}{l}=0.
$$
\end{proposition}

\begin{proof}
For every $\tilde k>0$, by Theorem \ref{mainthm} we have that for $l$ small enough
and $\Gamma \in \ks_m(\Omb)$ with $\hu(\Gamma) \le l$
$$
\int_\Om f(x,\nabla u)\,dx \le \int_\Om f(x,\nabla u_\Gamma)+\tilde{k} \hu(\Gamma).
$$
We deduce that 
$$
\limsup_{l \to 0^+}\frac{\ws(0)-\ws(l)}{l} \le \tilde{k}.
$$
Since $\tilde{k}$ is arbitrary, and since $\ws(l) \le \ws(0)$, we conclude that the result holds.
\end{proof}

\begin{remark}
\label{stronggriffith}
{\rm {\bf (The case of strong singularities)}
If the elastic solution $u$ has a strong singularity at $x \in \Om$, then by Remark
\ref{strongfalse} we have that
$$
\lim_{l \to 0^+}\frac{\ws(0)-\ws(l)}{l}=+\infty,
$$
so that the elastic configuration is not in equilibrium in the framework of Griffith's theory.
}
\end{remark}

\subsection{Crack initiation in variational theories of crack propagation} 
As explained at the beginning of the section, we consider now {\it irreversible} quasistatic crack evolutions governed by a {\it static equilibrium} condition and an {\it energy balance}. More precisely if $\psi(t)$ is a time dependent boundary displacement, and  $u(t),\Gamma(t)$ are the displacement and the crack at time $t$ relative to $\psi(t)$, we assume that the pair $(u(t),\Gamma(t))$ satisfies  the following properties:
\begin{itemize}
\item[$(a)$] {\it Irreversibility}: $\Gamma(t)$ is increasing in
time, i.e.,
$\Gamma(t_1) \subseteq \Gamma(t_2)$
for all $0 \le t_1 \le t_2 \le T$;
\vskip4pt
\item[$(b)$] {\it Static equilibrium}: if $t>0$, $\Es(u(t),\Gamma(t)) \le
\Es(u,H)$ for all cracks $H$ such that
$\cup_{s<t} \Gamma(s) \subseteq H$ and all displacements $v:\Om \setminus H \to \R$ with $v=\psi(t)$ on $\partial_D \Om \setminus H$;
\vskip4pt
\item[$(c)$] {\it Energy balance}: the total energy
$\Es(u(t),\Gamma(t))$
is absolutely continuous in time, and it satisfies
$$
\Es(u(t),\Gamma(t))=\Es(u(0),\Gamma(0))+\int_0^t \int_\Om 
\partial f(x,\nabla u(\tau)) \nabla \dot{\psi}(\tau)\,dx\,d\tau.
$$
\end{itemize}
\par
Condition $(a)$ stands for the {\it irreversibility} of
the evolution: the crack can only increase in time, i.e., no {\it healing processes} are admitted. Condition $(b)$ asserts that the pair $(u(t),\Gamma(t))$ is a {\it unilateral minimizer} of the total energy, i.e., it is a minimum among all configuration with larger cracks. In particular $u(t)$ is the elastic deformation relative to the boundary datum $\psi(t)$ in the domain $\Om \setminus \Gamma(t)$, i.e., $u(t)$ satisfies equation \eqref{elequation} with $\Gamma=\Gamma(t)$ and $\psi=\psi(t)$. Finally, notice that under suitable regularity assumptions on $u(t)$ and $\Gamma(t)$, condition $(c)$ can be rewritten as 
\begin{multline}
\label{balance}
\int_\Om f(x,\nabla u(t))\,dx-
\int_\Om f(x,\nabla u(s))\,dx+k\hu(K(t) \setminus K(s))\\
=\int_s^t \int_{\partial_D \Om \setminus K(\tau)} 
\frac{\partial}{\partial n} f(x,\nabla u(\tau)) \cdot \dot{\psi}(\tau)\,d\hu(x)\,d\tau.
\end{multline}
Therefore the {\it energy balance} condition states that the sum of the variation of the bulk energy and of the dissipation due to the creation of a new crack is equal to the work inserted in the system by the boundary datum $\psi$. We refer the reader to 
the paper by Francfort and Marigo \cite{FM} for further details on crack evolutions satisfying conditions $(a)$, $(b)$ and $(c)$ (see also Mielke \cite{Mi} for a connection with the theory of rate-independent processes).
\par
In order to treat the problem of crack initiation, we consider as in \cite{FM} the case in which $f$ is $p$-homogeneous in the gradient, i.e., for all $x \in \Om$, $\xi \in \R^2$ and $t>0$
$$
f(x,t\xi)=t^pf(x,\xi).
$$
\par
We consider a time dependent boundary displacement of the form $t \to t \psi$, where $t \in [0,T]$ and $\psi \in W^{1,p}(\Om) \cap L^\infty(\Om)$ is a given function. 
We refer the reader to the paper by Dal Maso and Toader \cite{DMT} for the existence of a quasistatic crack evolution satisfying $(a)$, $(b)$ and $(c)$ within this mathematical setting. 
\par
Since we are dealing with a crack initiation problem, and since $\psi(0)=0$, we assume that $(u(0),\Gamma(0))=(0,\emptyset)$. Notice that if $v$ denotes the elastic displacement associated to the boundary datum $\psi$, then $tv$ is the elastic displacement associated to $t \psi$. Then from the static equilibrium condition, comparing $(u(t),\Gamma(t))$ with $(tv,\Gamma(t))$, we have that for all $t \in ]0,T]$
\begin{equation}
\label{pivotineq2}
\int_\Om f(x,\nabla u(t))\,dx \le 
t^p\int_\Om f(x,\nabla v)\,dx.
\end{equation}
Finally, since we can replace $t\psi$ by $tu$ in the energy balance condition, 
we can write
\begin{equation}
\label{pivotineq}
\int_\Om f(x,\nabla u(t))\,dx+k\hu(\Gamma(t))=\int_0^t \int_\Om
\partial f(x,\nabla u(\tau))\cdot \nabla v\,dx\,d\tau.
\end{equation}
This implies that
\begin{equation}
\label{energybal}
\int_\Om f(x,\nabla u(t))\,dx+k\hu(\Gamma(t)) \le t^p\int_\Om f(x,\nabla v)\,dx.
\end{equation}
In fact, from the inequality
$$
xy \le f(x)+f^*(y),
$$ 
and taking into account the $p$-homogeneity of $f$, we obtain the following H\"older type inequality
\begin{equation}
\label{holdertype}
\int_\Omega zw \le p^{1/p}q^{1/q} \left( \int_\Omega f(x,z)\,dx \right)^{1/p}
\left( \int_\Omega f^*(x,w)\,dx \right)^{1/q},
\end{equation}
where $q=p':=p/(p-1)$. Taking $z=\nabla v$ and $w=\partial f(x,\nabla u(\tau))$ 
we get
\begin{multline*}
\int_\Om \partial f(x,\nabla u(\tau)) \nabla v\,dx 
\le
p^{1/p}q^{1/q}\left(\int_\Om f(x,\nabla v)\,dx\right)^{\frac{1}{p}} 
\left(\int_\Om f^*(\partial f(x,\nabla u(\tau)))\,dx\right)^{\frac{1}{q}} \\
=
p^{1/p}q^{1/q}\left(\int_\Om f(x,\nabla v)\,dx\right)^{\frac{1}{p}}  (p-1)^{\frac{1}{q}} 
\left(\int_\Om f(x,\nabla u(\tau))\,dx\right)^{\frac{1}{q}}. 
\end{multline*}
In view of \eqref{pivotineq2} we have
$$
\int_\Om \partial f(x,\nabla u(\tau)) \nabla v\,dx \le
\tau^{\frac{p}{q}} p^{1/p}q^{1/q}(p-1)^{\frac{1}{q}} \int_\Om f(x,\nabla v)\,dx.
$$
Integrating from $0$ to $t$ we obtain
$$
\int_0^t\int_\Om \partial f(x,\nabla u(\tau)) \nabla v\,dx 
\le
p^{1/p}q^{1/q}p^{-1} (p-1)^{\frac{1}{q}} t^p\int_\Om f(x,\nabla v)\,dx.
$$
Since
$$
p^{1/p}q^{1/q}p^{-1} (p-1)^{\frac{1}{q}}=1,
$$
by \eqref{pivotineq} we conclude that \eqref{energybal} holds.
\par
Notice that we can rescale \eqref{energybal} obtaining for $t$ small
\begin{equation}
\label{minimalityrescaled}
\int_\Om f(x,\nabla v(t))\,dx+\hu(\Gamma(t)) \le 
\int_\Om f(x,\nabla v)\,dx,
\end{equation}
where $v(t):=\frac{1}{t}u(t)$ is the displacement associated to $\Gamma(t)$ and $\psi$.
\par
As noticed by Francfort and Marigo in \cite{FM}, if $T$ is large enough, a crack will appear during the evolution, i.e., $\Gamma(s) \not=\emptyset$ for some $s \in ]0,T[$.
In fact if $T$ is such that 
\begin{equation}
\label{timebound}
\hu(\partial_D \Om)<T^p\int_\Om f(x,\nabla v)\,dx,
\end{equation}
then we get that creating a crack along $\partial_D \Om$ is more convenient that deforming $\Om$ elastically.
We are now in a position to state the first crack initiation result. 

\begin{theorem}
\label{weakinitiation}
Let us assume that $T$ satisfies \eqref{timebound}, and let us suppose that the elastic displacement $v$ associated to the boundary datum $\psi$ has uniformly weak singularities in $\Om$, i.e., it satisfies \eqref{weaksing}. Then the crack initiation is brutal, i.e., there exists a positive time $t_i\in ]0,T]$ such that $\Gamma(t)=\emptyset$ for every $t\le t_i$, and $\hs^1(\Gamma(t))>l^*$ for all $t \in ]t_i,T]$ for some $l^*>0$ depending only on $\Om$, $f$, $k$, $m$ and $\psi$.
\end{theorem}

\begin{proof}
Notice that by \eqref{energybal} we have $\hu(\Gamma(t)) \to 0$ as $t \to 0$. 
By Theorem \ref{mainthm}, we get that for $t$ small the elastic solution $v$ is energetically convenient with respect to $(v(t),\Gamma(t))$, where $v(t)$ is the displacement associated to $\Gamma(t)$ and $\psi$.
But this is against \eqref{minimalityrescaled} unless $\hu(\Gamma(t))=0$. In view of \eqref{timebound} we deduce that there exists  $t_i\in]0,T[$ such that $\Gamma(t)=\emptyset$ for every $t\le t_i$ and $\hs^1(\Gamma(t))>0$ for all $t \in ]t_i,T]$. 
\par
In order to prove that $\hs^1(\Gamma(t))>l^*$ for all $t \in ]t_i,T]$ and for some $l^*>0$,
notice that by \eqref{energybal} we deduce that
$$
\int_\Om f(x,\nabla v(t))\,dx+kT^{-p}\hu(\Gamma(t)) \le \int_\Om f(x,\nabla v)\,dx.
$$
Then by Theorem \ref{mainthm} we deduce that
$$
\liminf_{t \searrow t_i}\hu(\Gamma(t)) \ge l^*
$$ 
for some positive constant $l^*$ depending only on $\Om$, $f$, $kT^{-p}$, $m$ and $\psi$. The proof is thus concluded.
\end{proof}

Let us now study the crack initiation in the case in which the elastic displacement $u$ associated to the boundary datum $\psi$ has {\it strong singularities}. We recall that a point $x \in \Omb$ is a point of {\it strong singularity} for $v$ if 
$$
\limsup_{r \to 0} \frac{1}{r}\int_{B_r(x) \cap \Omb}|\nabla v|^p\,dx=+\infty.
$$
It is well expected that a during a loading process, a crack will appear at a point of strong singularity. The following theorem establishes this fact for a general quasistatic crack evolution satisfying properties $(a)$, $(b)$ and $(c)$.

\begin{theorem}
\label{stronginiziation}
Let us suppose that the elastic displacement $v$ associated to the boundary datum $\psi$ has a strong singularity at $x \in \Om$, and that $v$ has uniformly weak singularities in $\Omb \setminus B_r(x)$ for every $r>0$. Then we have that $\hu(\Gamma(t))>0$ for all $t \in ]0,T]$, and the crack starts at the point $\{x\}$, i.e.,
$$
x \in \bigcap_{t>0}\Gamma(t).
$$
Moreover the crack departs with zero speed, i.e.,
\begin{equation}
\label{speed}
\lim_{t \to 0}\frac{\hu(\Gamma(t))}{t}=0.
\end{equation}

\end{theorem}

\begin{proof}
In view of Remark \ref{strongfalse}, and since $x$ is a point of strong singularity for $v$, we have that for a positive time $t>0$ the elastic displacement $tv$ relative to the boundary datum $t\psi$ cannot satisfy condition $(b)$. As a consequence, we deduce that $\hu(\Gamma(t))>0$ for all $t \in ]0,T]$.
\par
Let us come to the properties of $\Gamma$ at time $t=0$. By \eqref{energybal} we have $$
\lim_{t \to 0}\hu(\Gamma(t))=0,
$$ 
so that
$$
\bigcap_{t>0}\Gamma(t)=\{y_1,\dots,y_h\},
$$
with $h \le m$. Let us suppose by contradiction that $y_i \not=x$ for all $i=1,\dots,h$. Then there exists $r>0$ such that for $t$ small enough we have $\Gamma(t) \subseteq \Omb \setminus B_r(x)$.
Since $v$ has uniformly weak singularities in $A:=\Om \setminus \bar B_r(x)$, and in view of Proposition \ref{mainthmrmk}, we have that inequality \eqref{minimalityrescaled} implies that $\hu(\Gamma(t))=0$ for $t$ small, which is a contradiction.
\par
Finally, in order to prove \eqref{speed}, we rescale \eqref{energybal} obtaining
$$
\int_\Om f(x,\nabla v(t)) \,dx+\frac{1}{t^p}\hu(\Gamma(t)) \le \int_\Om f(x,\nabla v)\,dx,
$$
where $v(t)$ is the displacement associated to $\Gamma(t)$ and $\psi$. We deduce  that $\hu(\Gamma(t)) \le C t^p$ for some constant $C$, and so \eqref{speed}
easily follows.
\end{proof}

\appendix
\section{The two dimensional $SBV$ case}
\label{SBVappendix}
The aim of this appendix is to prove a minimality result in the lines of Theorem \ref{mainthm} which does not require an a priori bound on the number of the connected components of the admissible cracks. Small cracks are still not energetically convenient if the gradient of the elastic solution of problem \eqref{elasticsol} and the related stress are continuous in $\Omb$. This condition excludes however that the elastic configuration presents (weak) singularities.
\par
In order to make the mathematical setting of this section precise, we need to recall some facts about {\it rectifiable sets} and the functional space $SBV$ of {\it special functions with bounded variation}. We refer the reader to \cite{AFP} for a complete treatment of these subjects.  
\par
A set $\Gamma \subseteq \R^N$ is rectifiable if there exists $N_0 \subseteq \Gamma$ with $\hn(N_0)=0$, and a sequence $(M_i)_{i \in \N}$ of $C^1$-submanifolds of $\R^N$ such that
$$
\Gamma \setminus N_0 \subseteq \bigcup_{i \in \N}M_i.
$$
For every $x \in \Gamma \cap M_i$, we define the normal to $\Gamma$ at $x$ as $n_{M_i}(x)$. It turns out that the normal is well defined (up to the sign) for $\hn$-a.e. $x \in \Gamma$.
\par
Let $U \subseteq \R^N$ be an open bounded set with Lipschitz boundary.
$SBV(U)$ is the set of functions $u \in L^1(U)$ such that the distributional derivative $Du$ is a Radon measure which, for every open set $A \subseteq U$, can be represented as
$$
Du(A)=\int_A \nabla u\,dx+\int_{A \cap S(u)}[u](x) \nu \,d\hn(x),
$$
where $\nabla u$ is the approximate differential of $u$, $S(u)$ is the set of jump
of $u$ (which is a rectifiable set), $\nu(x)$ is the normal to $S(u)$ at $x$, and $[u](x)$ is the jump of $u$ at $x$.
\par
For every $p \in ]1,+\infty[$ we set
$$
SBV^p(U):=\{u \in SBV(U)\,:\, \nabla u \in L^p(U,\R^N),\,\hn(S(u))<+\infty\}.
$$
If $u \in SBV(U)$, then $u$ admits a trace $\gamma(u)$ on $\partial U$ which is characterized by the relation (see \cite[Theorem 3.87]{AFP})
$$
\lim_{r \to 0}r^{-N}\int_{\Om \cap B_r(x)}|u(y)-\gamma(u)(x)|\,dy=0
\qquad
\text{for $\hn$-a.e. }x \in \partial U.
$$
We will denote the trace $\gamma(u)$ on $\partial U$ again by $u$. If $\Gamma \subseteq U$ is rectifiable and oriented by a normal vector field $n$, then we can define the traces $\gamma^+_\Gamma(u)$ and $\gamma^-_\Gamma(u)$ of $u \in SBV(U)$ on $\Gamma$ (see \cite[Theorem 3.77]{AFP}) which are characterized by the relations
$$
\lim_{r \to 0}r^{-N}\int_{\Om \cap B^\pm_r(x)}|u(y)-\gamma^\pm_\Gamma(u)(x)|\,dy=0
\qquad
\text{for $\hn$-a.e. }x \in \Gamma,
$$
where $B^\pm_r(x):=\{y \in B_r(x)\,:\, (y-x)\cdot n \gtrless 0\}$.
\par
A set $E \subseteq U$ has finite perimeter in $U$ if the characteristic function $1_E$ belongs to $SBV(U)$. We denote by $\partial^* E$ the set of jumps of $1_E$.
\par
Let us now come to our problem of local minimality. Let $\Om \subseteq \R^2$ be open, connected and with Lipschitz boundary. Let $\partial_D \Om \subseteq \partial \Om$ be open in the relative topology, and let $\partial_N \Om:=\partial \Om \setminus \partial \Om$. Let $f(x,\xi)$ be a Carath\'eodory function satisfying \eqref{fstruttura}, \eqref{fzero} and \eqref{fcrescita}. Let us moreover assume that the boundary displacement on $\partial_D \Om$ is given by the trace of a continuous function $\psi \in C^0(\Omb) \cap W^{1,p}(\Om)$. We denote by $u$ the elastic solution relative to $\psi$, namely the solution to the problem 
$$
\min \left\{ \int_\Om f(x,\nabla u)\,dx\,:\, u \in W^{1,p}(\Om), 
u=\psi \text{ on }\partial_D \Om\right\}.
$$
\par
The class of admissible cracks we consider is 
\begin{equation}
\label{crackgenerali}
\rs(\Omb):=\{\Gamma \subseteq \Omb\,:\, \Gamma \text{ is $1$-rectifiable and }\hu(\Gamma)<+\infty\}.
\end{equation}
\par
Let us come to the class of admissible displacements relative to a crack $\Gamma$
and to the boundary displacement $\psi$. Since $\Gamma$ is not supposed to be closed, the Sobolev space $W^{1,p}(\Om \setminus \Gamma)$ is not well defined. So we consider as class of admissible displacements the functions $u \in SBV^p(\Om)$ such that $\Sg{\psi}{u} \subseteq \Gamma$, where
$$
\Sg{\psi}{u}:=S(u) \cup \{x \in \partial_D \Om\,:\, u(x) \not= \psi(x)\},
$$
and the inequality on $\partial_D \Om$ is intended for the traces. Notice that if $\Gamma$ is closed  then $u \in W^{1,p}(\Om \setminus \Gamma)$, and $u=\psi$ on $\partial_D \Om \setminus \Gamma$.
\par
The displacement $u_\Gamma \in SBV^p(\Om) \cap L^\infty(\Om)$ associated to $\Gamma$ and $\psi$ is a solution of the minimum problem
\begin{equation}
\label{sbvpb}
\min \left\{ \int_\Om f(x,\nabla u)\,dx\,:\, u \in SBV^p(\Om), \Sg{\psi}{u} \subseteq \Gamma
\right\}.
\end{equation}
The proof of the existence of $u_\Gamma$ is standard: it relies on Ambrosio's compactness and lower semicontinuity theorem \cite{A2}, together with a truncation argument.
\par
The main result of the section is the following one.

\begin{theorem}
\label{sbvresult}
Let $u$ be the elastic displacement relative to $\psi \in W^{1,p}(\Om) \cap C^0(\Omb)$, and let us assume that $u$ satisfies
\begin{equation}
\label{continuityassumptions}
\nabla u  \in C^0(\Omb;\R^2)
\qquad\text{and}\qquad
\sigma:=\partial_\xi f(x,\nabla u) \in C^0(\Omb;\R^2).
\end{equation}
Then there exists a critical length $l^*>0$ depending only on $\Om$, $f$ and $\psi$ such that for all $\Gamma \in \rs(\Omb)$ with $\hu(\Gamma)<l^*$ we have
$$
\int_\Om f(x,\nabla u)\,dx< \int_\Om f(x,\nabla u_{\Gamma})\,dx+k\hu(\Gamma).
$$
\end{theorem}

In order to prove Proposition \ref{sbvresult} we need the following lemma.

\begin{lemma}
\label{sbvestimatelem}
For every $\Gamma \in \rs(\Omb)$ we have
\begin{equation}
\label{sbvestimate}
\int_\Om [f(x,\nabla u)-f(x,\nabla u_\Gamma)]\,dx \le \int_\Gamma \sigma\cdot n (u_\Gamma^+-u_\Gamma^-)\,d\hu,
\end{equation}
where $\sigma$ is the stress of the elastic displacement $u$ defined in \eqref{continuityassumptions}, $u_{\Gamma}$ is a minimum of \eqref{sbvpb}, and $u^\pm_\Gamma$ are the traces of $u_\Gamma$ on $\Gamma$.
\end{lemma}

\begin{proof}
By the convexity of $f$ we have
\begin{equation}
\label{convexineq}
\int_\Om [f(x,\nabla u)-f(x,\nabla u_\Gamma)]\,dx \le 
\int_\Om \partial f(x,\nabla u)(\nabla u-\nabla u_\Gamma)\,dx=\int_\Om \sigma (\nabla u-\nabla u_\Gamma)\,dx.
\end{equation}
We can assume that $\psi$ is defined on $\R^2$, i.e., $\psi \in W^{1,p}(\R^2) \cap C^0(\R^2)$. Let $B$ be a ball centred at $0$ such that $\Omb \subseteq B$.
Let us set $\Om':=B \setminus \partial_N \Om$. We can extend $u$ and $u_\Gamma$ to $\Om'$ setting $u=u_\Gamma=\psi$ on $B \setminus \Omb$. Let us consider $u_n \in C^1(\Om')$ with $u_n=\psi$ on $B \setminus \Omb$ and such that
$$
Du_n \to Du \quad\text{ strictly in the sense of measures,}
$$
that is (see \cite[Theorem 3.9]{AFP})
$$
\lim_{n \to +\infty}\int_{\Om'} \varphi dDu_n=\int_{\Om'} \varphi dDu
\qquad\text{ for all }\varphi \in C^0(\Om') \cap L^{\infty}(\Om').
$$
Since $u$ is a minimum for problem \eqref{sbvpb} with $\Gamma=\emptyset$, and $u_n=\psi$ on $\partial_D \Om$, then we have that $u-u_n$ is an admissible variation for $u$ so that for all $n \in \N$
\begin{equation}
\label{euler}
0=\int_\Om \sigma (\nabla u-\nabla u_n)\,dx
=\int_\Om \sigma \,dD(u-u_n).
\end{equation}
Since by assumption $\sigma \in C^0(\Omb)$, we can extend $\sigma$ to $\Om'$ in such a way that $\sigma \in C^0(\Om') \cap L^{\infty}(\Om')$. Notice that since $u=u_n=\psi$ on $B \setminus \Omb$, from \eqref{euler} we have
$$
\int_{\Om'}\sigma\,dD(u_n-u)=0.
$$
Then by strict convergence we deduce
\begin{multline}
\label{appineq1}
0=\lim_{n \to \infty} \int_{\Om'} \sigma \,dD(u-u_n)
= \int_{\Om'} \sigma \nabla u\,dx-\int_{\Om'}\sigma \,dDu_\Gamma\\
=\int_{\Om'}\sigma (\nabla u-\nabla u_\Gamma)\,dx-
\int_{S(u_\Gamma)} \sigma \cdot n (u_\Gamma^+-u_\Gamma^-)\,d\hu.
\end{multline}
Since $S(u_\Gamma) \subseteq \Gamma$, $u_\Gamma^+=u_\Gamma^-$ on $\Gamma \setminus S(u_\Gamma)$, and $\nabla u=\nabla u_\Gamma=\nabla \psi$ on $B \setminus \Omb$, by \eqref{appineq1} we deduce
$$
\int_{\Gamma} \sigma \cdot n (u_\Gamma^+-u_\Gamma^-)\,d\hu=
\int_{\Om'}\sigma (\nabla u-\nabla u_\Gamma)\,dx=\int_{\Om}\sigma (\nabla u-\nabla u_\Gamma)\,dx,
$$
so that, in view of \eqref{convexineq}, we have that \eqref{sbvestimate} follows.
\end{proof}

We are now in a position to prove the minimality result in the $SBV$ context.

\begin{proof}[Proof of Proposition \ref{sbvresult}]
Let us consider $l>0$, and let $\Gamma \in \rs(\Omb)$ be a minimum for the functional
$$
\fs(\Gamma):=\int_\Om f(x,\nabla u_{\Gamma})\,dx+k\hu(\Gamma)
$$
among the cracks $\Gamma \in \rs(\Omb)$ such that $\hu(\Gamma) \le l$. The existence of such a $\Gamma$ follows minimizing the functional
$$
F(u):=\int_\Om f(x,\nabla u_{\Gamma})\,dx+k\hu(\Sg{\psi}{u})
$$
among all $u \in SBV(\Om)$ with $\hu(\Sg{\psi}{u}) \le l$, and choosing $\Gamma:=S(u)$.
Proposition \ref{sbvresult} will be proved if we show that $\Gamma=\emptyset$ for $l$ small enough.
\par
Let us establish some qualitative properties of $\Gamma$ that are useful in the proof. Let $E \subseteq \Om$ be a set with finite perimeter such that $\partial^* E \subseteq \Gamma$. We can assume that
\begin{equation}
\label{uGammazero}
u_{\Gamma}=0
\qquad
\text{on }E.
\end{equation}
In fact we have that 
$$
\tilde u_{\Gamma}:=
\begin{cases}
u_{\Gamma}      &\text{in }\Om \setminus E\\
0      &\text{in }E
\end{cases}
$$
belongs to $SBV(\Om)$, and it is an admissible displacement for $\Gamma$ and $\psi$ (see \cite[Theorem 3.84]{AFP}) with
$$
\int_\Om f(x,\nabla \tilde u_{\Gamma})\,dx \le
\int_\Om f(x,\nabla u_{\Gamma})\,dx.
$$
We conclude that $\tilde u_{\Gamma}$ is a minimum energy displacement relative to $\Gamma$ and $\psi$, so that \eqref{uGammazero} follows.
\par
Let $E_1$ and $E_2$ be two sets in $\Om$ with finite perimeter such that
$\partial^* E_1 \subseteq \Gamma$ and $\partial^* E_2 \subseteq \Gamma$.
We have that 
\begin{equation}
\label{h1zero}
\hu(\partial^* E_1 \cap \partial^* E_2)=0.
\end{equation}
In fact, if by contradiction $\hu(\partial^* E_1 \cap \partial^* E_2)>0$,  let us consider the crack
$$
\tilde \Gamma:= \Gamma \setminus (\partial^* E_1 \cap \partial^* E_2).
$$
Notice that $S(u_\Gamma) \subseteq \tilde \Gamma$. In fact $S(u_{\Gamma}) \subseteq \Gamma$, and $u_{\Gamma}=0$ on $E_1 \cup E_2$ in view of \eqref{uGammazero} so that 
$$
\hu\left( S(u_{\Gamma}) \cap (\partial^* E_1 \cap \partial^* E_2) \right)=0.
$$ 
Since $\hu(\tilde \Gamma)<\hu(\Gamma)\le l$ and 
$u_{\Gamma}$ is an admissible displacement for $\tilde \Gamma$,
we get that $\fs(\tilde\Gamma)<\fs(\Gamma)$, and this is against the minimality of $\Gamma$. Then \eqref{h1zero} holds.
\par
As a consequence of \eqref{h1zero}, there exists a family $(E_i)_{i \in \N}$ of sets with finite perimeter in $\Om$ such that $\partial^* E_i \subseteq \Gamma$, $\hu(\partial^* E_h \cap \partial^* E_k)=0$ for every $h \not=k$, and such that if $E$ is a set with finite perimeter in $\Om$ such that $\partial^* E \subseteq \Gamma$, then $E= E_h$ for some $h\in\N$. Let us set
\begin{equation}
\label{defgammastar}
\Gamma^*:=\bigcup_{i=0}^\infty \partial^* E_i,
\end{equation}
and let us assume that $n$ denotes the outward normal to $E_i$. Since the stress $\sigma$ is a divergence free vector field, in view of the Generalized Gauss-Green formula for sets with finite perimeter (see \cite[Theorem 3.36]{AFP}) we have for all $i \in \N$
\begin{equation}
\label{divfreeEk}
\int_{\partial^* E_i} \sigma \cdot n\,d\hu=\int_{E_i} {\rm div}\sigma\,dx=0.
\end{equation}
\par
In order to prove our minimality result, by contradiction let us assume that there exists $l_h \to 0$ and $\Gamma_h$ such that setting $u_h:=u_{\Gamma_h}$ we have
\begin{equation}
\label{contradiction}
\int_\Om f(x,\nabla u_h)\,dx+k\hu(\Gamma_h)<\int_\Om f(x,\nabla u)\,dx.
\end{equation}
By Ambrosio's lower semicontinuity Theorem \cite{A2} and by \eqref{contradiction},
we deduce that for every open set $A \subseteq \Om$
\begin{equation}
\label{convenergia}
\lim_{h \to +\infty} \int_A f(x,\nabla u_h)\,dx=\int_A f(x,\nabla u)\,dx.
\end{equation}
Let $(E^h_i)_{i \in \N}$ and $\Gamma^*_h$ be the sets associated to $\Gamma_h$ described above. We claim that
\begin{equation}
\label{jumpest}
|[u_h]|:=|u_h^+-u_h^-| \to 0
\qquad
\text{uniformly on $\Gamma_h \setminus \Gamma_h^*$ as }h \to +\infty
\end{equation}
and
\begin{equation}
\label{essoscest}
\text{ess-sup}_{\partial^* E^h_i}(u_h^+)-\text{ess-inf}_{\partial^* E^h_i}(u_h^+) \to 0
\qquad
\text{uniformly in $i$ as }h \to +\infty.
\end{equation}
In view of \eqref{jumpest} and \eqref{essoscest}, the proof of the proposition is readily concluded. In fact, given $\eps>0$, and choosing $h$ so large that 
$$
|[u_h]| \le \eps
\qquad
\text{ on } \Gamma_h \setminus \Gamma_h^*
$$
and for all $i \in \N$
$$
\text{ess-sup}_{\partial^* E^h_i}(u_h^+)-\text{ess-inf}_{\partial^* E^h_i}(u_h^+) \le \eps, 
$$
by Lemma \ref{sbvestimatelem} and in view also of \eqref{uGammazero} and of \eqref{divfreeEk} we have that 
\begin{align*}
\int_\Om &[f(x,\nabla u)-f(x,\nabla u_h)]\,dx \le \int_{\Gamma_h} \sigma\cdot n (u_h^+-u_h^-)\,d\hu\\
&\le \int_{\Gamma_h \setminus \Gamma_h^*} \sigma\cdot n (u_h^+-u_h^-)\,d\hu+
\sum_{i=0}^\infty \int_{\partial^* E^h_i} \sigma\cdot n u_h^+\,d\hu \\
&\le \eps \|\sigma\|_{\infty} \hu\left( \Gamma_h \setminus \Gamma_h^* \right)+
\sum_{i=0}^\infty (\text{ess-inf}_{\partial^*E^h_i}u_h^+)\int_{\partial^* E^h_i} \sigma\cdot n \,d\hu + \eps \|\sigma\|_{\infty} \sum_{i=0}^\infty \hu(\partial^*E^h_i) \\
&= \eps  \|\sigma\|_{\infty} \left(\hu\left( \Gamma_h \setminus \Gamma_h^* \right)+\hu\left(\Gamma_h^* \right)\right)= \eps  \|\sigma\|_{\infty} \hu\left( \Gamma_h\right),
\end{align*}
and this is against \eqref{contradiction}.
\par
In order to conclude the proof, we have to prove the claims \eqref{jumpest} and \eqref{essoscest}. Let us consider \eqref{jumpest}, the proof of the other claim being similar.
Let us assume that there exists $\delta>0$ and $x_h \in \Gamma_h \setminus \Gamma_h^*$ with 
\begin{equation}
\label{claim1}
|[u_h](x_h)| \ge \delta>0. 
\end{equation}
Up to a subsequence we have $x_h \to \bar x \in \Omb$. Let us assume that $\bar x \in \Om$. For $h$ large enough, and for $r$ small we have $\bar B_r(x_h) \subseteq \Om$. Notice that for every $r$ such that $\partial B_{r}(x_h) \cap \Gamma_h=\emptyset$,
we have that $u_h \in W^{1,p}(\partial B_r(x_h))$. Moreover by the maximum principle we have that 
\begin{equation}
\label{maxprinc}
\max_{\partial B_{r}(x_h)}u_h-\min_{\partial B_{r}(x_h)}u_h=M_h-m_h>\delta.
\end{equation}
In fact otherwise, we can consider $\tilde u_h$ defined as
$$
\tilde u_h:=
\begin{cases}
u_h      &\text{outside }B_r(x_h),\\
\max\{\min\{u_h,M_h\}, m_h\}     &\text{inside }B_r(x_h).
\end{cases}
$$
Since $x_h \in \Gamma_h \setminus \Gamma_h^*$, and in view of \eqref{claim1} we deduce that 
$$
|\{u_h \not= \tilde u_h\} \setminus \bigcup_{i \in \N} E^h_i|>0
$$
so that 
$$
\int_\Om f(x,\nabla \tilde u_h)\,dx<\int_\Om f(x,\nabla u_h)\,dx
$$
which is against the minimality of $u_h$. Then \eqref{maxprinc} holds.
\par
Let $\partial_{\vartheta} u_h$ denote the angular derivative of $u_h$, i.e., $\partial_{\vartheta} u_h:=\frac{d}{d\vartheta}u_h(x_h+r\cos \vartheta,y_h+r\sin\vartheta)$. Setting $C_h:=\{s \in [0,r]\,:\, \hu \left( \partial B_s(x_h) \cap \Gamma_h \right)=0\}$, by \eqref{maxprinc} we have for every $r \in C_h$ 
$$
\int_0^{2\pi}|\partial_\vartheta u_h|^p\,d\vartheta \ge 2\pi\delta^p.
$$
Since $|C_h| \ge r-\hu(\Gamma_h)$, we deduce that
\begin{multline}
\label{tetaintegrale}
\int_{B_r(x_h)}|\nabla u_h|^p\,dx \ge \int_0^r \int_0^{2\pi}s^{1-p}|\partial_\vartheta u_h|^p \,d\vartheta\,ds \\
\ge \int_{C_h}\int_0^{2\pi} s^{1-p}|\partial_\vartheta u_h|^p \,d\vartheta\,ds \ge 2\pi \delta^p \int_{C_h}s^{1-p}\,ds \ge 2\pi\delta^p \int_{\hu(\Gamma_h)}^r s^{1-p}\,ds.
\end{multline}
Let us distinguish two cases, namely $p \ge 2$ and $1<p<2$. If $p \ge 2$, 
choosing $r=2\hu(\Gamma_h)$ we obtain
$$
\liminf_{h \to +\infty} \int_{B_{r_h}(x_h)}f(x,\nabla u_h)\,dx \ge 2\pi\delta^p \ln2.
$$
But this is against \eqref{convenergia}: in fact for all $r$ such that $B_r(\bar x) \subseteq \Om$ by \eqref{convenergia} we have
$$
\int_{B_r(\bar x)}f(x,\nabla u)\,dx=
\lim_{h \to +\infty}
\int_{B_r(\bar x)}f(x,\nabla u_h)\,dx \ge 
\lim_{h \to +\infty}
\int_{B_{r_h}(x_h)}f(x,\nabla u_h)\,dx \ge 2\pi\delta^p \ln2,
$$
and this gives a contradiction for $r$ small enough.
If $1<p<2$, then we have
\begin{equation*}
\liminf_{h \to +\infty}  \int_{B_{r}(x_h)}f(x,\nabla u_h)\,dx \ge
\liminf_{h \to +\infty }\frac{2\pi\delta^p}{2-p} \left(r^{2-p}-\hu(\Gamma_h)^{2-p}\right) = \frac{2\pi\delta^p}{2-p} r^{2-p}
\end{equation*}
from which by \eqref{convenergia} we deduce that
\begin{equation}
\label{sing}
\int_{B_r(\bar x)}f(x,\nabla u)\,dx \ge Cr^{2-p}
\end{equation}
for some $C>0$. Since $\nabla u \in C^0(\Omb;\R^2)$ and $f$ satisfies \eqref{fcrescita}
we get that
$$
\int_{B_r(\bar x)}f(x,\nabla u)\,dx \le \tilde C r^2
$$
for some $\tilde C>0$, which together with \eqref{sing} gives a contradiction.
\par
The case in which $\bar x \in \partial \Om$ can be treated almost in the same way as
the case $\bar x \in \Om$. In fact it is sufficient to choose $r$ so small that
$$
\max_{B_r(x_h) \cap \partial_D \Om}\psi
-\min_{B_r(x_h) \cap \partial_D \Om}\psi<\delta,
$$
and to take into account the fact that there exists a constant $C$ depending only on $\Om$ such that
$$
\hu(\partial B_{r_h}(x_h) \cap \Om) \ge C r_h.
$$
In this way, integrations involved in \eqref{tetaintegrale} can be performed on a set of angle $\vartheta$ which has a positive measure uniformly bounded from below, and the contradiction follows by the same arguments used above.
\end{proof}

\def\bu{\mathbf{u}}%

\section{Uniform Poincar\'e and Poincar\'e-Korn inequalities}
\label{app:PK}

In this section, we show that very basic arguments lead to the uniform
Poincar\'e and Poincar\'e-Korn inequalities that are needed respectively
in Sections~\ref{mainresultsec} and~\ref{planarelasticitysec}.
\par
To simplify, we will only
consider the case of $N$-dimensional domains of the
form
\[
Q_f=
\left\{ x=(x',x_N)\in \R^N\,:\, 0\leq x_i\leq 1\,,i=1,\dots,N-1\,,
0 \le x_N\leq f(x') \right\}
\]
where $f:Q'\to [1,M]$ is a $L$-Lipschitz function. Here 
$Q'$ is the $(N-1)$-dimensional cube $(0,1)^{N-1}$ and
$L>0$ and $M>1$ are fixed constants.
The adaption of the argument
that we will present here to the ``real'' cases that are useful in the paper
is straightforward. Also, for simplicity, we consider here the
``linear'' case $p=2$. However, the proofs would be identical with any
other exponent $p\in (1,+\infty)$.
\par
With a slight abuse in the notation we also identify $Q'$ with
the base of $Q_f$, that is, the subset
$(0,1)^{N-1}\times \{0\}$ of $\partial Q_f$.
We show that the following result holds:
\begin{proposition}\label{th:PK}
There exists a constant $C>0$ depending only on $L$ and $M$
such that
\begin{itemize}
\item[(i)] For any $u \in H^1(Q_f)$ with $u=0$ on $Q'$,
$\|u\|_{L^2(Q_f)}\leq C\|\nabla u\|_{L^2(Q_f)}$;
\item[(ii)]  For any $u \in H^1(Q_f)$ with $\int_{Q_f} u(x)\,dx=0$,
$\|u\|_{L^2(Q_f)}\leq C\|\nabla u\|_{L^2(Q_f)}$;
\item[(iii)] For any $\bu \in H^1(Q_f;\R^N)$ with $\bu=0$ on $Q'$,
$\|\bu\|_{L^2(Q_f)}\leq C\|e(\bu)\|_{L^2(Q_f)}$;
\item[(iv)]  For any $\bu \in H^1(Q_f;\R^N)$ with both
 $\int_{Q_f} \bu(x)\,dx=0$ and $\int_{Q_f} x\times \bu(x)\,dx=0$,
one has
$\|\bu\|_{L^2(Q_f)}\leq C\|e(\bu)\|_{L^2(Q_f)}$.
\end{itemize}
\end{proposition}
In the last assertion, $x\times \bu$ is the skew-symmetric matrix
$(x_i \bu_j-x_j\bu_i)_{i,j=1}^N$, and the condition means that
$\bu$ is orthogonal (in $L^2$) to the rigid motions (of the form $\mathbf{a}
+Bx$ with $B$ skew-symmetric).
\par
Let us sketch the proof of this proposition. First of all, the proof
of point (i) is standard (by integration along vertical lines starting
from $Q'$) and it is well known that the constant $C$, in this case,
only depends on $M$ ($f$ could then be any l.s.c.\ function below $M$).
In the same way, the proof of (iii) is significantly simpler than
the proof of (iv) (note however that it does require that $f$ is
Lipschitz and $C$ will depend on both $M$ and $L$), and we will not discuss it.
(See~\cite{BC02} for a detailed proof, in dimension two).
\par
To prove (ii) one first establishes the following inequality:
there exists $C_0$ depending only on $M$ such that for any
$f$ and any $u\in H^1(Q_f)$, one has
\begin{equation}\label{poin1}
\int_{Q_f} u(x)^2\,dx \ \leq\ C_0\left(
\int_A u(x)^2\,dx \,+\, \int_{Q_f} |\nabla u(x)|^2\right)
\,dx\,,
\end{equation}
where $A$ denotes the set $Q'\times (0,1)$ (the important fact here being
that $A$ is an open set that belongs to \textit{all} the domains $Q_f$,
for all admissible $f$).
The proof of~\eqref{poin1}, again, is standard. It relies on integration
along vertical lines starting from the base $Q'$ and on the obvious fact
that for any $b\in [1,M]$ and any $v\in C^1(0,b)$,
\begin{equation}\label{poin2}
\int_0^y v(t)^2\,dt \ \leq\ 2M\int_0^1 v(t)^2\,dt \,+\,2M^2\int_0^b v'(t)^2
\,dt\,.
\end{equation}
\par
Now, if (ii) is not true, it means that there exists
functions $f_n$ and $u_n$
with $1\leq f_n\leq M$, $f_n$ $L$-Lipschitz, $u_n\in H^1(Q_{f_n})$, $\int_{Q_{f_n}}u_n\,dx=0$ and
\[
\|u_n\|_{L^2(Q_{f_n})} \ \geq \ n\|\nabla u_n \|_{L^2(Q_{f_n})}
\]
for any $n$.
\par
Without loss of generality we may renormalize $u_n$ so that
$\int_A u_n^2\,dx=1$. Then, by~\eqref{poin1} we find
\[
\|\nabla u_n \|_{L^2(Q_{f_n})}\,\leq\,\frac{1}{n}
\|u_n\|_{L^2(Q_{f_n})}\,+\, \leq \frac{\sqrt{C_0}}{n}
\left( 1+\|\nabla u_n\|_{L^2(Q_{f_n})}\right)\,.
\]
If we extend both $u_n$ and $\nabla u_n$ with the value
$0$ outside of $Q_{f_n}$, this inequality shows that
$\nabla u_n$ goes to zero strongly in $L^2(Q_M;\R^N)$ ($Q_M=
Q\times (0,M)$), and, up to a subsequence, that there exists
$u\in L^2(Q_M)$ such that $u_n\rightharpoonup u$ weakly in $L^2(Q_M)$.
\par
On the other hand, by Rellich's theorem, $u_n\to u$ strongly
in $L^2(A)$ and since $\nabla u=0$ in $A$ and $\int_A u^2=1$,
$u$ is the constant $\pm1/\sqrt{|A|}$.
\par
Now, since the functions $f_n$ are uniformly equibounded and equicontinuous,
up to a further subsequence, we may assume also that $f_n$ converges to
some $f$ uniformly. It is now easy to check that $u\in H^1(Q_f)$, $u=0$
outside of $Q_f$, and $\nabla u=0$ (the limit of $\nabla u_n$)
in $Q_f$ so that $u$ is a constant in $Q_f$.
We deduce that $u=(1/\sqrt{|A|})\chi_{Q_f}$. Now, for each $n$, one
had $\int_{Q_{f_n}}u_n\,dx=\int_{Q_M}u_n\,dx=0$, hence in the limit
$\int_{Q_M}u\,dx=0\neq |Q_f|/\sqrt{|A|}$, a contradiction. Hence (ii)
must be true.
\par
We observe here that (ii) holds in fact as long as $f$ belongs to a
fixed set of functions which is compact in $C^0(Q',[1,+\infty))$
(the constant $C$ depending only on this compact set). The case
of $L$-Lipschitz functions uniformly bounded by the constant $M$ is
a particular case. On the other hand, for the Poincar\'e-Korn inequalities
(iii) and (iv), the fact that the functions $f$ are uniformly
Lipschitz seems
to be essential, as we now show.
\par
Let us now prove (iv). It is enough to show that the vectorial version
of~\eqref{poin1} holds, that is,
\begin{equation}\label{poink1}
\int_{Q_f} |\bu(x)|^2\,dx \ \leq\ C_0\left(
\int_A |\bu(x)|^2\,dx \,+\, \int_{Q_f} |e(\bu)(x)|^2\right)
\,dx\,.
\end{equation}
This will be shown, again, by integration along lines and using~\eqref{poin2},
however, this time, it is not sufficient to consider only vertical lines
starting from $Q'$. Indeed, one has for any smooth vectorial field
$\bu\in C^1(Q_f)$
that 
$$
d(\bu(x+s\xi)\cdot \xi)/ds= (e(\bu)(x+s\xi)\xi)\cdot \xi
$$
for any $x\in Q_f$, $\xi\in \mathbb{S}^{N-1}$, and $s\in \R$ such that
$x+s\xi\in Q_f$. Hence, integration along vertical lines will control
the component $u_N$ of $\bu=(u_1,\dots,u_N)$. To control the other
components, one needs to integrate along lines in at least
$N-1$ other independent directions (as is done in~\cite{BC02}).
\par
Given $p>\max\{L,2M\}$ let us consider, for any $i=1,\dots,N-1$,
the vectors
\[
\xi_i^\pm =\frac{1}{\sqrt{1+p^2}}(0,\dots,0,\pm 1,0,\dots, p)
\] where $\pm 1$
appears at the $i$th position.
Given $i\leq N-1$, if we considers the lines starting from $Q'$
in the direction $\xi_i^+$, we see that they ``see'' all points
$x\in Q_f$ with $x_i\geq 1/2$. On the other hand, the lines starting
from $Q'$ in the direction $\xi_i^-$ ``see'' all the points
with $x_i\leq 1/2$. Integrating along these lines and using~\eqref{poin2},
one controls the $L^2$-norms on one half of the domain $Q_f$ of
$(u_i+pu_N)/\sqrt{1+p^2}$ and, on the other half, of
$(-u_i+pu_N)/\sqrt{1+p^2}$. Together with the control of
$\int_{Q_f} u_N^2\,dx$ obtained previously, this shows that one
can control $\int_{Q_f} u_i^2\,dx$ with the right-hand side of~\eqref{poink1}.
Repeating this argument for all $i$, we find that~\eqref{poink1}
holds, with now a constant that depends on $M$ and $L$, through $p$.
\par
We deduce, exactly as before, that~(iv) holds.\qed

%% \vskip20pt\noindent{\bf Acknowledgments.}
\section*{Acknowledgments}
The authors wish to thank Gilles Francfort, Jean-Jacques
Marigo and Vincenzo Nesi for interesting discussions. The first and third
authors are partially funded by the MULTIMAT European network
MRTN-CT\_2004-505226.

\end{document}